\documentclass[11pt,reqno]{amsart} 
 \usepackage{amsmath,amssymb,amsthm}
 \usepackage{graphicx}
 \usepackage{cite}
 \usepackage[mathscr]{eucal}

 \theoremstyle{plain}
 \newtheorem{theorem}{Theorem}[section]  
 \newtheorem{corollary}[theorem]{Corollary}  
 \newtheorem{lemma}[theorem]{Lemma}  
 
 \newtheorem*{theorem*}{Theorem}

 \theoremstyle{plain}  
 
 \newtheorem{definition}[theorem]{Definition}

 \newtheoremstyle{citing}
   {3pt}
   {3pt}
   {\itshape}
   {}
   {\bfseries}
   {.}
   {.5em}
   {\thmnote{#3}}

 \theoremstyle{citing}

\numberwithin{equation}{section}
\allowdisplaybreaks[1]

\newlength{\intwidth}
\makeatletter
\DeclareRobustCommand{\cpvint}[2]
    {\mathop{%
       \text{%
         \settowidth{\intwidth}{%
           \ifx\ilimits@\displaylimits
             $\int_{#1}^{#2}$%
           \else
             $\int$%
           \fi}%
         \makebox[0pt][l]{\makebox[\intwidth]{$\text{C}$}}%
         $\int_{#1}^{#2}$}}}
\makeatother

\makeatletter
\DeclareRobustCommand{\cpvintsmall}[2]
    {\mathop{%
       \text{%
         \settowidth{\intwidth}{%
           \ifx\ilimits@\displaylimits
             $\int_{#1}^{#2}$%
           \else
             $\int$%
           \fi}%
         \makebox[0pt][l]{\makebox[\intwidth]{$\text{{\tiny C}}$}}%
         $\int_{#1}^{#2}$}}}
\makeatother

\newcommand{\D }{\Delta }

\newcommand{\xpfeil}{\xrightarrow}  
  
\newcommand{\rand}{\partial} 
\newcommand{\where}{:\:}  
\newcommand{\iso}{\cong}  
  
\newcommand{\sgn}{\text{sgn}}

\newcommand{\cupl}{\mathop{\cup}\limits}  
\newcommand{\capl}{\mathop{\cap}\limits}

\newcommand{\sumg}{\mathop{\sum}\limits}  
  
\newcommand{\aequi}{\Longleftrightarrow}

\newcommand{\nz}{{\mathbb N}}

\newcommand{\rz}{{\mathbb R}}  
\newcommand{\zz}{{\mathbb Z}}  
 
\newcommand{\eps}{\varepsilon}  
\renewcommand{\phi}{\varphi} 
\newcommand{\eval}{\vert}
\newcommand{\cross}{\times}

\begin{document}
 
\title[Closed magnetic geodesics] {Closed Magnetic Geodesics on $S^2$}
\author{Matthias Schneider}
\address{Ruprecht-Karls-Universit\"at\\
         Im Neuenheimer Feld 288\\
         69120 Heidelberg, Germany\\}
\email{mschneid@mathi.uni-heidelberg.de} 
\date{April 8, 2010}  

\begin{abstract}
We give existence results for simple closed
curves with prescribed geodesic curvature on $S^{2}$, which correspond to periodic orbits
of a charge in a magnetic field.
\end{abstract}

\maketitle

\section{Introduction}
\label{sec:introduction}
The trajectory of a charged particle on an 
orientable Riemannian surface $(N,g)$ in a magnetic field
given by the magnetic field form $\Omega = k\, dA$, 
where $k:N \to \rz$ is the magnitude of the magnetic field and $dA$ is the area form on $N$,
corresponds to a curve $\gamma$ on $N$ that solves
\begin{align}
\label{eq:magnetic_geodesic}
D_{t,g} \dot \gamma = k(\gamma) J_{g}(\gamma) \dot \gamma
\end{align}
where $D_{t,g}$ is the covariant derivative with respect to $g$,
and $J_g(x)$ is the rotation by $\pi/2$ in $T_xN$ measured with $g$ and the orientation chosen 
on $N$. A curve $\gamma$ in $N$ that solves (\ref{eq:magnetic_geodesic}) will be called a {\em ($k$-)magnetic geodesic}.  
We refer to \cite{MR890489,MR1432462,MR2036336} for the Hamiltonian description
of the motion of a charge in a magnetic field.
Taking the scalar product of \eqref{eq:magnetic_geodesic} with $\dot\gamma$ we see that
if $\gamma$ is a magnetic geodesic then 
$(\gamma,\dot\gamma)$ lies on the energy level $E_c:=\{(x,V)\in TN \where |V|_g=c\}$.\\
The geodesic curvature $k_g(\gamma,t)$ of 
an immersed curve $\gamma$ at $t$ is defined by
\begin{align*}
k_g(\gamma,t) := |\dot\gamma(t)|_{g}^{-2}\big\langle \big(D_{t,g} \dot\gamma\big)(t) , N_{g}(\gamma(t))\big\rangle_{g},  
\end{align*}
where $N_{g}(\gamma(t))$ denotes
the unit normal of $\gamma$ at $t$ given by
\begin{align*}
N_{g}(\gamma(t)) := |\dot\gamma(t)|_{g}^{-1} J_{g}(\gamma(t))\dot\gamma(t).  
\end{align*}
By \eqref{eq:magnetic_geodesic}, a nonconstant curve $\gamma$ on $E_c$ is 
a $k$-magnetic geodesic if and only if
its geodesic curvature $k_g(\gamma,t)$ is given by $k(\gamma(t))/c$.
We will take advantage of the latter description and consider
the equation
\begin{align}
\label{eq:1}
D_{t,g} \dot \gamma = |\dot \gamma|_{g} k(\gamma) J_{g}(\gamma) \dot \gamma.  
\end{align}
We call equation (\ref{eq:1}) the {\em prescribed geodesic curvature equation}, as its solutions $\gamma$ 
are constant speed curves with geodesic curvature $k_g(\gamma,t)$ given by $k(\gamma(t))$.
For fixed $k$ and $c>0$ the equations \eqref{eq:magnetic_geodesic} and \eqref{eq:1} are equivalent in the following
sense: If $\gamma$ is a nonconstant solution of \eqref{eq:1} with $k$ replaced by
$k/c$, then the curve 
$t \mapsto \gamma(ct/|\dot \gamma|_g)$
is a $k$-magnetic geodesic on $E_c$ and a $k$-magnetic geodesic on $E_c$
solves \eqref{eq:1} with $k$ replaced by $k/c$. We emphasize that unless $k\equiv 0$
the solutions of \eqref{eq:magnetic_geodesic} lying in different $E_c$ are
not reparametrizations of each other.\\
We study the existence of closed curves with prescribed geodesic curvature 
or equivalently the existence of periodic magnetic geodesics on prescribed energy levels 
$E_c$.\\
There are different approaches to this problem, the Morse-Novikov theory
for (possibly multi-valued) variational functionals (see \cite{MR1185286,MR730159,MR1133303}), the theory of dynamical
systems using methods
from symplectic geometry (see \cite{MR2324797,MR2060021,MR902290,MR890489,MR1432462,MR1417851,MR2208800}) 
and Aubry-Mather's theory (see \cite{MR2036336}).\\
We suggest a new approach, instead of looking for critical
points of the (possibly multivalued) action functional
we consider solutions to \eqref{eq:1} as zeros of 
the vector field $X_{k,g}$ defined on the Sobolev space $H^{2,2}(S^1,N)$ as follows:
For $\gamma \in H^{2,2}(S^{1},N)$ we let $X_{k,g}(\gamma)$ be the unique weak solution of   
\begin{align}
\label{eq:def_vector_field}
\big(-D_{t,g}^{2} + 1\big)X_{k,g}(\gamma)= 
-D_{t,g} \dot\gamma + |\dot \gamma|_{g} k(\gamma)J_{g}(\gamma)\dot\gamma   
\end{align}
in $T_\gamma H^{2,2}(S^{1},N)$.  
The uniqueness implies that any zero of $X_{k,g}$ is a weak solution of \eqref{eq:1}
which is a classical solution in $C^{2}(S^{1},N)$ applying standard regularity
theory. 
The vector field $X_{k,g}$ as well as the set of solutions to \eqref{eq:1}
is invariant under a circle action:
For $\theta \in S^{1}=\rz/\zz$ and $\gamma \in H^{2,2}(S^{1},N)$
we define $\theta*\gamma \in H^{2,2}(S^{1},N)$
by 
\begin{align*}
\theta*\gamma(t) = \gamma(t+\theta).  
\end{align*}
Moreover, for $V \in T_\gamma H^{2,2}(S^{1},N)$ we let
\begin{align*}
\theta*V := V(\cdot+\theta) \in T_{\theta*\gamma} H^{2,2}(S^{1},N).  
\end{align*}
Then $X_{k,g}(\theta * \gamma) = \theta*X_{k,g}(\gamma)$ for any $\gamma \in H^{2,2}(S^{1},N)$
and $\theta\in S^{1}$. Thus, any zero gives rise to a $S^{1}$-orbit of zeros and we say
that two solutions $\gamma_1$ and $\gamma_2$ of (\ref{eq:1}) are (geometrically) distinct, if
$S^{1}*\gamma_1 \neq S^{1}*\gamma_2$.\\
We will apply this approach to the case $N=S^2$, equipped with a smooth metric $g$, and
$k$ a positive smooth function on $S^{2}$. 
We shall prove
\begin{theorem}
\label{thm_existence}
Let $g$ be a smooth metric and $k$ a positive smooth function on $S^2$.
Suppose that one of the following three assumptions is satisfied, 
\begin{align}
\label{eq:cond_inj}
4\inf(k) 
\ge \big(inj(g)\big)^{-1}\big(2\pi+(\sup K_g^{-})vol(S^2,g)\big),\\
\label{eq:cond_K_pos}
K_g>0 \text{ and } 2\inf(k)
 \ge \sup(K_g)^{\frac12},\\
\label{eq:cond_K_pinch} 
\sup(K_g) <4\inf(K_g),
\end{align}
where $K_g$ denotes the Gauss curvature, $K_g^-:= -\min(K_g,0)$, and $inj(g)$ the injectivity radius
of $(S^2,g)$. 
Then there are at least two simple solutions of \eqref{eq:1}
in $C^2(S^1,S^2)$.
\end{theorem}
Concerning the existence of closed $k$-magnetic geodesics for a positive smooth function $k$ 
on $(S^{2},g)$ the following is known (see \cite{MR902290,MR1432462})
\begin{itemize}
\item[(i)]
if $c>0$ is sufficiently small, then $E_c$ contains
two simple closed magnetic geodesics,
\item[(ii)] if $g$ is sufficiently close to the round metric $g_0$ and $k$ is sufficiently close to
a positive constant, then there is a closed magnetic geodesic
in every energy level $E_c$,
\item[(iii)]
if $c> 0$ is sufficiently large, then $E_c$ contains a closed magnetic geodesic.
\end{itemize}
Using the equivalence between (\ref{eq:magnetic_geodesic}) and (\ref{eq:1}) we obtain
from Theorem \ref{thm_existence}
\begin{corollary}
\label{cor_existence}
Let $g$ be a smooth metric, $k$ a positive smooth function on $S^2$,
and $c>0$.
Suppose that one of the following three assumptions is satisfied,
\begin{align}
\label{eq:cond_inj_c}
c \le 4 \big(\inf(k)\big)\big(inj(g)\big)\big(2\pi+(\sup K_g^{-})vol(S^2,g)\big)^{-1},\\
\label{eq:cond_K_pos_c}
K_g>0 \text{ and } c\le 2\inf(k) \big(\sup(K_g)\big)^{-\frac12},\\ 
\notag
\sup(K_g) <4\inf(K_g).
\end{align}
Then $E_c$ contains at least two simple closed magnetic geodesics.
\end{corollary}
Condition \eqref{eq:cond_inj_c} should be compared to the existence results in $(i)$ and
gives bounds on the required smallness of $c$ in terms of geometric quantities.
To show that 
\eqref{eq:cond_inj_c} is useful despite the implicit definition of $inj(g)$, 
we apply an estimate of $inj(g)$ in \cite{MR1330918} and obtain
\eqref{eq:cond_K_pos_c} as a special case.
The pinching condition \eqref{eq:cond_K_pinch} extends the existence result in $(ii)$
and shows for instance that on the round sphere there are two simple closed curves of prescribed
geodesic curvature $k$ for any positive function $k$, which gives a partial
solution to a problem posed by Arnold in \cite[1994-35,1996-18]{MR2078115} 
concerning the existence of closed magnetic geodesics 
on $S^{2}$ on every energy level $E_c$.\\
By the famous Lusternik-Schnirelmann theorem there are at least three simple closed
geodesics on every Riemannian two sphere $(S^{2},g)$. As a by-product of our analysis we show that in general,
even if $k$ is very close to $0$, there are no more than two simple
closed magnetic geodesics on $S^{2}$ in $E_1$ (see also \cite[Sec. 7]{MR2060021}).
\begin{theorem}
\label{thm_example}
Let $g_0$ be the round metric on $S^{2}$. For any positive constant $k_0>0$ 
there is a smooth function $k$ on $S^{2}$, which can be chosen arbitrarily close to $k_0$, such that
there are exactly two simple solutions of (\ref{eq:1}).  
\end{theorem}
The proof of our existence results is organized as follows.
After setting up notation in Section \ref{s:preliminiaries} and introducing the classes of maps
and spaces needed for our analysis 
we define in Section \ref{sec:s1_euler}
a $S^1$-equivariant Poincar\'{e}-Hopf index or $S^1$-degree, $$\chi_{S^1}(X,M) \in \zz,$$ where $M$ is a $S^1$-invariant
subset of prime curves in $H^{2,2}(S^1,S^2)$ and $X$ belongs to
a class of $S^1$-invariant vector fields. 
The index $\chi_{S^1}(X,M)$ 
is related to the extension
of the Leray-Schauder degree to intrinsic nonlinear problems in \cite{MR493919,MR0277009}
and is used 
combined with the apriori estimates in Section \ref{sec:apriori} to count simple periodic solutions of (\ref{eq:1}).
We remark that the standard degree $\chi(X,M)$, 
that does not take the $S^1$ invariance into account, 
vanishes as it detects only fixed points under the $S^1$-action, i.e. constant solutions.
Equivariant degree theories have been defined and applied to
differential equations by many authors, we refer to
\cite{MR2036376,MR0209600,MR817033,MR1984999,MR2282341} and the references therein.
However, we do not see how to apply these results directly to (\ref{eq:1}).\\
The vector field $X_{k,g}$ corresponding
to the prescribed geodesic curvature problem falls into the class of vector fields,
where our $S^1$-degree is defined. In Section \ref{sec:deg_poincare} we show
that the $S^1$-degree of an isolated zero orbit of $X_{k,g}$
is given by $-i(P,\theta)$, where $i(P,\theta)$ 
denotes the fixed point index of the Poincar\'{e} map of the corresponding magnetic geodesic.
Since the Poincar\'{e} map is area preserving we obtain from \cite{MR0353372,MR0346846}
that the $S^1$-degree of an isolated zero orbit is bounded below by $-1$.\\     
Section \ref{sec:unperturbed}
is devoted to the computation
of $\chi_{S^1}(X_{k_0,g_0},M)$,
where $k_0$ is a positive constant, $g_0$ is the round metric of $S^2$, and
$M$ is the set of simple regular curves in $H^{2,2}(S^1,S^2)$.
We call equation \eqref{eq:1} with $k\equiv k_0$ and $g=g_0$ the unperturbed problem,
which is analyzed in detail. The set of simple solution
to the unperturbed problem is given by circles of latitude
of radius $(1+k_0^2)^{-1/2}$ with respect to an arbitrarily chosen north pole.
To compute the $S^1$-degree we slightly perturb 
the constant function $k_0$ and end up with exactly two nondegenerate
solutions of degree $-1$. This implies that
$\chi_{S^1}(X_{k_0,g_0},M)=-2$.\\  
Section \ref{sec:apriori} contains the apriori estimates showing
that the set of simple solutions to \eqref{eq:1} is compact in $M$
under each of the assumptions \eqref{eq:cond_inj}-\eqref{eq:cond_K_pinch}.
This yields together with the perturbative analysis in Section \ref{sec:unperturbed}
the proof of Theorem \ref{thm_example} and allows to construct an admissible homotopy
of vector fields between $X_{k_0,g_0}$ and $X_{k,g}$ whenever $k$ and $g$
satisfy the assumptions of Theorem \ref{thm_existence}.   
The homotopy invariance of the $S^{1}$-equivariant Poincar\'{e}-Hopf index then shows
\begin{align*}
\chi_{S^1}(X_{k,g},M)=\chi_{S^1}(X_{k_0,g_0},M)=-2.  
\end{align*}
Since the $S^1$-degree of an isolated zero orbit is always larger than $-1$
there are at least two simple solutions of \eqref{eq:1}.   
The existence result is given in Section \ref{sec:existence}.

\subsection*{Acknowledgements}
I would like to thank Professor Friedrich Tomi for valuable discussions
and the unknown referees for several suggestions which helped to improve the article, 
in particular for drawing my attention to \cite{MR0353372,MR0346846}.

\section{Preliminaries}
\label{s:preliminiaries}
Let $S^{2}=\rand B_1(0) \subset \rz^{3}$ be the standard round sphere with induced metric $g_0$
and orientation such that the rotation $J_{g_0}(y)$ is given for $y \in S^2$ by
\begin{align*}
J_{g_0}(y)(v) := y \cross v \text{ for all }v \in T_y S^{2},   
\end{align*}
where $\cross$ denotes the cross product in $\rz^{3}$. 
If we equip $S^{2}$ with a general Riemannian metric $g$, then
the rotation by $\pi/2$ measured with $g$ is given by 
\begin{align*}
J_{g}(y) v = \big(G(y)\big)^{-1} J_{g_0}(y)\big(G(y)\big)v  \quad \forall v \in T_y S^{2},
\end{align*}
where $G(y)$ denotes 
a positive symmetric map $G(y) \in \mathcal{L}(T_y S^{2})$ satisfying
\begin{align*}
\langle v, w\rangle_{T_yS^2,g} =
\langle G(y)v,G(y)w\rangle_{T_yS^2,g_0} \quad \forall v,w \in T_y S^{2}.  
\end{align*}
We consider for $m\in \nz_0$ the set of Sobolev functions
\begin{align*}
H^{m,2}(S^{1},S^{2}) := \{\gamma \in H^{m,2}(S^{1},\rz^{3}) \where \gamma(t) \in \rand B_1(0) 
\text{ for a.e. } t \in S^{1}.\}   
\end{align*}
For $m \ge 1$ the set
$H^{m,2}(S^{1},S^{2})$ is a sub-manifold of the Hilbert space $H^{m,2}(S^{1},\rz^{3})$
and is contained in $C^{m-1}(S^{1},\rz^{3})$.
Hence, if $m\ge 1$ then $\gamma \in H^{m,2}(S^{1}, S^{2})$ 
satisfies 
$\gamma(t) \in \rand B_1(0)$ for all $t\in S^{1}$. 
In this case the tangent space $T_\gamma H^{m,2}(S^{1}, S^{2})$ 
of $H^{m,2}(S^{1}, S^{2})$ 
at $\gamma \in H^{m,2}(S^{1}, S^{2})$ is given by
\begin{align*}
T_\gamma H^{m,2}(S^{1}, S^{2}) := 
\{V \in H^{m,2}(S^{1},\rz^{3})\where V(t) \in T_{\gamma(t)}S^{2} \text{ for all }t \in S^{1}\}.  
\end{align*}
For $m=0$ the set $H^{0,2}(S^{1},S^{2})=L^{2}(S^{1},S^{2})$ fails to be a manifold. 
In this case we define for $\gamma \in H^{1,2}(S^{1},S^{2})$ the space 
$T_\gamma L^{2}(S^{1},S^{2})$ by
\begin{align*}
T_\gamma L^{2}(S^{1},S^{2}):= 
\{V \in L^{2}(S^{1},\rz^{3})\where V(t) \in T_{\gamma(t)}S^{2} \text{ for a.e. }t \in S^{1}\}.  
\end{align*}
A metric $g$ on $S^{2}$ induces a metric on $H^{m,2}(S^{1},S^{2})$ for $m\ge 1$ 
by setting for $\gamma \in H^{m,2}(S^{1},S^{2})$
and $V,\,W \in T_\gamma H^{m,2}(S^{1}, S^{2})$
\begin{align*}
\langle W,V\rangle_{T_\gamma H^{m,2}(S^{1}, S^{2}),g} 
:= \int_{S^{1}}\Big\langle &\big((-1)^{\llcorner \frac{m}{2}\lrcorner}(D_{t,g})^{m}+1\big)V(t),\\
&\big((-1)^{\llcorner \frac{m}{2}\lrcorner}(D_{t,g})^{m}+1\big)W(t)\Big\rangle_{\gamma(t),g}\, dt,   
\end{align*}
where $\llcorner m/2\lrcorner$ denotes the largest integer that does not exceed $m/2$.\\
Let $X$ be a differentiable vector field on $H^{2,2}(S^1,S^2)$. Then 
the covariant (Frechet) derivative $D_{g} X$,
\begin{align*}
D_{g} X: \: TH^{2,2}(S^{1},S^{2})) \to TH^{2,2}(S^{1},S^{2}),   
\end{align*}
of the vector field $X$ with respect to the metric induced by $g$ is defined as follows:
For $\gamma \in H^{2,2}(S^{1},S^{2})$ and $V \in T_\gamma  H^{2,2}(S^{1},S^{2})$ we consider a $C^{1}$-curve
\begin{align*}
(-\eps,\eps) \ni s \mapsto \gamma_s \in H^{2,2}(S^{1},S^{2})  
\end{align*}
satisfying 
\begin{align*}
\gamma_0 = \gamma \text{ and } \frac{d}{ds}\gamma_s\eval_{s=0} = V,  
\end{align*}
and define
\begin{align*}
D_{g} X\eval_\gamma[V](t) := D_{s,g}\Big(X\big(\gamma_s(t)\big)\Big)\eval_{s=0}.  
\end{align*}
For the vector field theory on infinite dimensional manifolds it is convenient to work
with Rothe maps instead of compact perturbations of the identity, because the class
of Rothe maps is open in the space of linear continuous maps.
We recall the definition and properties of Rothe maps given in \cite{MR493919} for the sake of the readers convenience. 
For a Banach space $E$ we denote by $\mathcal{GL}(E)$ the set of invertible maps in $\mathcal{L}(E)$
and by $\mathcal{S}(E)$ the set
\begin{align*}
\mathcal{S}(E)=\{T \in \mathcal{GL}(E) \where (tT+(1-t)I)\in \mathcal{GL}(E) \text{ for all }t \in [0,1]\}.  
\end{align*}
Then the set of Rothe maps $\mathcal{R}(E)$ is defined by
\begin{align*}
\mathcal{R}(E) := \{A \in \mathcal{L}(E)\where A=T+C,\, T\in \mathcal{S}(E) \text{ and } C \text{ compact}\}.  
\end{align*}
The set $\mathcal{R}(E)$ is open in $\mathcal{L}(E)$ and consists of Fredholm operators of index $0$. 
Moreover, $\mathcal{GR}(E):= \mathcal{R}(E) \cap \mathcal{GL}(E)$ 
has two components, $\mathcal{GR}^{\pm}(E)$, with $I \in \mathcal{GR}^{+}(E)$.
For $A \in \mathcal{GR}(E)$ we let
\begin{align*}
\sgn A =
\begin{cases}
+1 &\text{if }A \in \mathcal{GR}^{+}(E),\\
-1 &\text{if }A \in \mathcal{GR}^{-}(E).
\end{cases}
\end{align*}
If $A=I+C \in \mathcal{GL}(E)$, where $C$ is compact, then $A \in \mathcal{GR}(E)$ and $\sgn A$ is given by the
the usual Leray-Schauder degree of $A$.\\
Since $g$ and $k$ are smooth, $X_{k,g}$ is a smooth vector field (see
\cite[Sec. 6]{MR0464304}) on the set $H^{2,2}_{reg}(S^{1},S^{2})$ of regular curves,
\begin{align*}
H^{2,2}_{reg}(S^{1},S^{2}) := \{\gamma \in H^{2,2}(S^{1},S^{2})\where \dot\gamma(t) \neq 0 \text{ for all }t\in S^{1}\}. 
\end{align*}
To compute $D_gX_{k,g}\eval_\gamma(V)$ we observe
\begin{align*}
D_{s,g}&\big(-D_{t,g}^2+1)X_{k,g}(\gamma_s)\big) =
D_{s,g}\big(-D_{t,g}\dot\gamma_s +|\dot\gamma_s|_g k(\gamma_s)J_{g}(\gamma_s)\dot\gamma_s\big)\\
&=-D_{t,g}^2 \frac{d\gamma_s}{ds} - R_g\big(\frac{d\gamma_s}{ds},\dot\gamma_s\big)\dot\gamma_s   
+D_{s,g}\big(|\dot\gamma_s|_g k(\gamma_s)J_{g}(\gamma_s)\dot\gamma_s\big).
\end{align*}
Evaluating at $s=0$ we obtain
\begin{align}
\label{eq:5}
D_{s,g}&\big(-D_{t,g}^2+1)X_{k,g}(\gamma_s)\big)\eval_{s=0} \notag \\
 &=
-D_{t,g}^2 V - R_g\big(V,\dot\gamma \big)\dot\gamma
+ |\dot\gamma|_g^{-1}\langle D_{t,g}V,\dot\gamma\rangle_g k(\gamma)J_{g}(\gamma)\dot\gamma
 \notag \\
&\quad +|\dot\gamma|_g \big(k'(\gamma)V\big)J_{g}(\gamma)\dot\gamma 
+|\dot\gamma|_g k(\gamma) \big(D_gJ_{g}\eval_{\gamma}V\big)\dot\gamma \notag\\
&\quad+|\dot\gamma|_g k(\gamma)J_{g}(\gamma) D_{t,g}V.
\end{align}
Moreover, we have
\begin{align}
\label{eq:6}
D_{s,g}&\big(-D_{t,g}^2+1)X_{k,g}(\gamma_s)\big)\eval_{s=0} \notag\\
&= -D_{s,g}D_{t,g}^2 X_{k,g}(\gamma_s)\eval_{s=0} +D_{s,g}X_{k,g}(\gamma_s)\eval_{s=0} \notag \\
&= \big(-D_{t,g}^2+1)D_{g} X_{k,g}\eval_\gamma(V)
-D_{t,g}\Big(R_g\big(V,\dot\gamma \big)X_{k,g}(\gamma)\Big)\notag \\
&\quad - R_g\big(V,\dot\gamma\big)D_{t,g}X_{k,g}(\gamma).  
\end{align}
Equating (\ref{eq:5}) and (\ref{eq:6}) at a critical point $\gamma$ of $X_{k,g}$ leads to
\begin{align}
\label{eq:dg_x_g_formula}
\big(-&D_{t,g}^2+1\big)D_{g} X_{k,g}\eval_\gamma(V) \notag\\
&= -D_{t,g}^2 V - R_g\big(V,\dot\gamma \big)\dot\gamma
+ |\dot\gamma|_g^{-1}\langle D_{t,g}V,\dot\gamma\rangle_g k(\gamma)J_{g}(\gamma)\dot\gamma
 \notag \\
&\quad +|\dot\gamma|_g \big(k'(\gamma)V\big)J_{g}(\gamma)\dot\gamma
+ |\dot\gamma|_g k(\gamma) \big(D_gJ_{g}\eval_{\gamma}V\big)\dot\gamma \notag\\
&\quad+|\dot\gamma|_g k(\gamma)J_{g}(\gamma) D_{t,g}V. 
\end{align}
We note that (see also \cite[Thm. 6.1]{MR493919})
\begin{align*}
\big(-D_{t,g}^2+1\big)D_{g} X_{k,g}\eval_\gamma(V) = (-D_{t,g}^2 +1)V +T(V),   
\end{align*}
where $T$ is a linear map from $T_\gamma H^{2,2}(S^{1},S^{2})$ to $T_\gamma L^{2}(S^{1},S^{2})$ that
depends only on $V$ and its first derivatives and is therefore compact.
Taking the inverse $(-D_{t,g}^2+1)^{-1}$ we deduce that $D_{g} X_{k,g}\eval_\gamma$
is the form $identity + compact$ and thus a Rothe map.\\
For $m\ge 1$ the exponential map $Exp_{g}: TH^{m,2}(S^{1},S^{2}) \to H^{m,2}(S^{1},S^{2})$
is defined for $\gamma \in H^{m,2}(S^{1},S^{2})$ and $V \in T_\gamma H^{m,2}(S^{1},S^{2})$ by
\begin{align*}
Exp_{\gamma,g}(V)(t) := Exp_{\gamma(t),g}(V(t)),  
\end{align*}
where $Exp_{z,g}$ denotes the exponential map on $(S^{2},g)$ at $z \in S^{2}$. Due to its pointwise definition
\begin{align*}
\theta*Exp_{\gamma,g}(V)(t) = Exp_{\theta*\gamma,g}(\theta*V)(t).    
\end{align*}

\section{The $S^{1}$-Poincar\'{e}-Hopf index}
\label{sec:s1_euler}
For $\gamma \in H^{2,2}(S^{1},S^{2})$ we define the form $\omega_{g}(\gamma) \in (T_\gamma H^{2,2}(S^{1},S^{2}))^{*}$ by
\begin{align*}
\omega_{g}(\gamma)(V) :&= \int_0^{1}\langle \dot \gamma(t), \big(-(D_{t,g})^2 + 1\big) V(t)\rangle_{\gamma(t),g}\,dt\\
&= \langle\dot\gamma, V\rangle_{T_\gamma H^{1,2}(S^{1},S^{2}),g}. 
\end{align*}
Approximating $\dot\gamma$ by vector fields contained in $T_\gamma
H^{2,2}(S^{1},S^{2})$, it is easy to see that $\omega_{g}(\gamma)\neq 0$, if $\gamma \neq {\rm const}$.
If $\gamma \in H^{3,2}(S^{1},S^{2})$, then $\omega_{g}(\gamma)$ extends to a linear form in
$(T_\gamma L^{2}(S^{1},S^{2}))^{*}$ by
\begin{align*}
\omega_g(\gamma)(V) := \langle \big(-(D_{t,g})^2 + 1\big)\dot\gamma, V\rangle_{T_\gamma L^{2}(S^{1},S^{2}),g}.
\end{align*}
From Riesz' representation theorem there is
$W_{g}(\gamma)\in T_\gamma H^{2,2}(S^{1},S^{2})$ such that
\begin{align*}
\omega_{g}(\gamma)(V) = \langle V, W_{g}(\gamma)\rangle_{T_\gamma H^{2,2}(S^{1},S^{2}),g}
\; \forall  V \in T_\gamma H^{2,2}(S^{1},S^{2}),  
\end{align*}
and
\begin{align}
\label{eq:w_g_gamma_h12}
\langle W_g(\gamma)\rangle^\perp = \langle \dot\gamma \rangle^{\perp,H^{1,2}}\cap T_\gamma H^{2,2}(S^1,S^2).  
\end{align}
Hence $$W_{g}(\gamma) = (-(D_{t,g})^2 + 1)^{-1} \dot \gamma$$ 
and $W_{g}$ is a $C^{2}$ vector field on $H^{2,2}(S^{1},S^{2})$.\\
The form $\omega_{g}(\gamma)$ and the vector $W_{g}(\gamma)$ are equivariant under the $S^{1}$-action in the sense that
for all $\theta \in S^{1}$ and $V\in T_\gamma H^{2,2}(S^{1},S^{2})$ we have
\begin{align*}
w_{\theta*\gamma,g}(\theta*V)=\omega_{g}(\gamma)(V) \text{ and } W_{\theta*\gamma,g}=\theta*W_{g}(\gamma).   
\end{align*}
Using the vector field $W_g$ we define a vector bundle $SH^{2,2}(S^{1},S^{2})$ by
\begin{align*}
SH^{2,2}(S^{1},S^{2}) := \{(\gamma,V)\in TH^{2,2}(S^{1},S^{2})
\where \gamma \neq {\rm const},\,V\in \langle W_g(\gamma)\rangle^\perp\}.   
\end{align*}
Note that $SH^{2,2}(S^{1},S^{2})$ is $S^{1}$-invariant, as
\begin{align*}
(\gamma,V) \in SH^{2,2}(S^{1},S^{2}) \implies (\theta*\gamma,\theta*V) \in SH^{2,2}(S^{1},S^{2})\; \forall \theta
\in S^{1}.   
\end{align*}
For $\gamma \in H^{2,2}(S^{1},S^{2})\setminus\{{\rm const}\}$ we consider the map 
\begin{align*}
\psi_{\gamma,g}:T_\gamma H^{2,2}(S^{1},S^{2}) \times T_\gamma H^{2,2}(S^{1},S^{2}) \to SH^{2,2}(S^{1},S^{2})    
\end{align*}
defined by
\begin{align}
\label{eq:def_psi}
\psi_{\gamma,g}(V,U) &:=
\Big(Exp_{\gamma,g}V,Proj_{\langle W_g(Exp_{\gamma,g}V)\rangle^{\perp}}\big(D Exp_{\gamma,g}\eval_V U\big)\Big).
\end{align}
The differential of $\psi_{\gamma,g}$ at $(0,0)$ is given by
\begin{align*}
D\psi_{\gamma,g}\eval_{(0,0)}(V,U) = 
(V,U-\|W_g(\gamma)\|^{-2} \langle U,W_g(\gamma)\rangle_{T_{\gamma} H^{2,2}(S^{1},S^{2}),g}W_g(\gamma)).  
\end{align*}
Consequently, there is $\delta=\delta(\gamma,g)>0$ such that 
$\psi_{\gamma,g}$ restricted to 
\begin{align*}
B_{\delta}(0)\times B_{\delta}(0)\cap \langle W_g(\gamma)\rangle^\perp \subset
T_{\gamma} H^{2,2}(S^{1},S^{2})\times T_{\gamma} H^{2,2}(S^{1},S^{2})  
\end{align*}
is a chart for the manifold $SH^{2,2}(S^{1},S^{2})$ at $(\gamma,0)$.
The construction is $S^{1}$-equivariant, for
\begin{align*}
\psi_{\theta*\gamma,g}(\theta*V,\theta*U)=\theta*\psi_{\gamma,g}(V,U) \; \forall \theta \in S^{1}
\end{align*}
and we may choose $\delta(\gamma,g)=\delta(\theta*\gamma,g)$ for all $\theta \in S^{1}$.
Shrinking $\delta(\gamma,g)$ we may assume, as $Exp_{\gamma,g}$ is also a chart 
for $H^{k,2}(S^1,S^2)$ with $1\le k\le 4$ and by \eqref{eq:w_g_gamma_h12},
\begin{align}
\label{eq:def_trans_1}
T_{Exp_{\gamma,g}(V)}H^{1,2}(S^{1},S^{2}) 
= \langle D_tExp_{\gamma,g}(V)\rangle \oplus DExp_{\gamma,g}\eval_V(\langle \dot\gamma\rangle^{\perp,H^{1,2}}),\\
\label{eq:def_trans_2}
T_{Exp_{\gamma,g}(V)}H^{2,2}(S^{1},S^{2}) 
= \langle W_g(Exp_{\gamma,g}(V))\rangle \oplus DExp_{\gamma,g}\eval_V(\langle W_g(\gamma)\rangle^{\perp}),\\
\label{eq:def_trans_3}
\text{Proj}_{\langle W_{g}(Exp_{\gamma,g}(V)\rangle^\perp}\circ D Exp_{\gamma,g}\eval_{V}:\,
\langle W_{g}(\gamma)\rangle^\perp \xpfeil{\iso}{} \langle W_{g}(Exp_{\gamma,g}(V)\rangle^\perp,
\end{align}
and the norm of the projections corresponding to the decompositions in (\ref{eq:def_trans_1}) and (\ref{eq:def_trans_2})
as well as the norm of the map in (\ref{eq:def_trans_3}) 
and its inverse are uniformly bounded with respect to $V$.\\
The $S^1$-action is only continuous but not differentiable on $H^{2,2}(S^1,S^2)$
as for instance the candidate for the differential of the map $\theta \to
\theta*\gamma$ at $\theta =0$, $\dot \gamma$, is in general only in $T_\gamma H^{1,2}(S^1,S^2)$.
We prove the existence of a slice of the $S^{1}$-action (see \cite[Lem. 2.2.8]{MR0478069} and the references
therein) at a curve $\gamma$ with higher regularity and obtain additional differentiability  
of the slice map. 
\begin{lemma}[Slice lemma]
\label{l:slice_lemma}
Let $\gamma \in H^{3,2}(S^{1},S^{2})$ be a prime curve, i.e.
a curve with trivial isotropy group $\{\theta \in S^1\where \theta*\gamma=\gamma\}$. 
Then there is an open neighborhood $\mathcal{U}$
of $0$ in $T_\gamma H^{2,2}(S^{1},S^{2})$, such that 
the map $$\Sigma_{\gamma,g}: S^{1} \times \mathcal{U}\cap \langle W_g(\gamma)\rangle^{\perp} \to
H^{2,2}(S^{1},S^{2}),$$ 
defined by
\begin{align*}
\Sigma_{\gamma,g}(\theta,V) := \theta*Exp_{\gamma,g}(V),  
\end{align*}
is a homeomorphism onto its range, which is open in $H^{2,2}(S^{1},S^{2})$. Moreover, the inverse
$(\Sigma_{\gamma,g})^{-1}$ satisfies
\begin{align*}
\text{Proj}_{S^{1}}\circ (\Sigma_{\gamma,g})^{-1}\in C^{2}\Big(\Sigma_{\gamma,g}\big(S^{1} \times \mathcal{U}\cap \langle
W_g(\gamma)\rangle^{\perp}\big),S^{1}\Big).  
\end{align*}
\end{lemma}
\begin{proof}
Fix a prime curve $\gamma \in H^{3,2}(S^{1},S^{2})$. 
We consider for $\delta_0>0$ the map
$$F_{\gamma,g}:B_{\delta_0}(0)\times B_{\delta_0}(0)\subset \rz/\zz \times T_\gamma H^{2,2}(S^{1},S^{2}) \to
\rz$$
defined by 
\begin{align*}
F_{\gamma,g}(\theta,V) := \omega_{g}(\gamma)\Big(Exp_{\gamma,g}^{-1}\big(\theta*Exp_{\gamma,g}(V)\big)\Big).
\end{align*}
Note that, as $S^{1}$ acts continuously on $H^{2,2}(S^{1},S^{2})$ and $Exp_{\gamma,g}$ is a local diffeomorphism,
after shrinking $\delta_0>0$ the map $F_{\gamma,g}$ is well defined. 
$Exp_{\gamma,g}$ is a smooth map, such that
for fixed $\theta$ the map $V \mapsto F_{\gamma,g}(\theta,V)$ is also smooth.
Moreover, since
$Exp_{\gamma,g}(V)$ is in $H^{2,2}(S^1,S^2)$
and $D Exp_{\gamma,g}\eval_V$ maps $L^{2}$ vector fields along 
$\gamma$ into $L^{2}$ vector fields along $Exp_{\gamma,g}(V)$, the map
\begin{align*}
\theta \mapsto Exp_{\gamma,g}^{-1}\big(\theta*Exp_{\gamma,g}(V)\big)  
\end{align*}
is $C^2$ from $B_{\delta_0}(0)\subset \rz/\zz$ to $T_\gamma L^{2}(S^{1},S^{2})$, the space of
$L^{2}$ vector fields along $\gamma$.
For $\gamma \in H^{3,2}(S^{1},S^{2})$ the form $\omega_{g}(\gamma)$ is in $(T_\gamma L^{2}(S^{1},S^{2}))^{*}$. 
Thus, $\theta \mapsto F_{\gamma,g}(\theta,V)$ is $C^2$ as well as $F_{\gamma,g}$.
Fix $V \in T_\gamma H^{2,2}(S^{1},S^{2})$. 
Since
\begin{align*}
D_\theta F_{\gamma,g}\eval_{(0,0)} = \omega_{g}(\gamma)(\dot\gamma),
\neq 0.
\end{align*}
by the implicit function theorem and after shrinking $\delta_0>0$ we get a unique $C^{2}$-map
\begin{align*}
\sigma_{\gamma,g}:B_{\delta_0}(0)\subset T_\gamma H^{2,2}(S^{1},S^{2}) \to B_{\delta_0}(0)\subset
\rz/\zz  
\end{align*}
such that
\begin{align*}
F_{\gamma, g}(\sigma_{\gamma,g}(V),V) \equiv 0 \text{ in } B_{\delta_0}(0)\subset T_\gamma H^{2,2}(S^{1},S^{2}).  
\end{align*}
Hence, we may define locally around $\gamma$ 
\begin{align*}
V_{\gamma,g}(\alpha) :&= Exp_{\gamma,g}^{-1}(\sigma_{\gamma,g}(Exp_{\gamma,g}^{-1}(\alpha))*\alpha)
\in \langle W_g(\gamma)\rangle^\perp.
\end{align*}
Using the uniqueness of $\sigma_{\gamma,g}$ and the fact that $\gamma$ is prime 
it is standard to see that $\Sigma_{\gamma,g}$
is injective and that the inverse is given locally around 
$\theta_0*\gamma$ for fixed $\theta_0 \in S^1$ by
\begin{align*}
\Sigma_{\gamma,g}^{-1}= (\theta_0,0)+
(-\sigma_{\gamma,g}\circ Exp_{\gamma,g}^{-1}\circ (-\theta_0*),V_{\gamma,g}\circ (-\theta_0*)).
\end{align*}
This finishes the proof.
\end{proof}
We will compute the Poincar\'{e}-Hopf index for the following class of vector fields.
\begin{definition}
Let $M$ be an open 
$S^1$-invariant subset of prime curves in $H^{2,2}(S^{1},S^2)$.
A $C^2$ vector field $X$ on $M$ is called $(M,g,S^{1})$-admissible, if
\begin{itemize}
\item[(1)] $X$ is $S^1$-equivariant, i.e. $X(\theta*\gamma)=\theta*X(\gamma)$ for all
$(\theta,\gamma)\in S^1\times M$.
\item[(2)] $X$ is proper in $M$, i.e. the set $\{\gamma \in M \where X(\gamma)=0 \}$ is compact,
\item[(3)] $X$ is orthogonal to $W_{g}$, i.e. $w_{g}(\gamma)(X(\gamma))=0$ for all $\gamma \in M$.
\item[(4)] $X$ is a Rothe field, i.e. 
if $X(S^{1}*\gamma)=0$ then 
\begin{align*}
D_{g}X\eval_\gamma &\in \mathscr{R}(T_\gamma H^{2,2}(S^1,S^2)) \text{ and }
\text{Proj}_{\langle W_g(\gamma)\rangle^\perp} \circ D_{g}X\eval_\gamma
\in \mathscr{R}(\langle W_g(\gamma)\rangle^\perp),
\end{align*}
\item[(5)] $X$ is elliptic, i.e. there is $\eps>0$ such that  
for all finite sets of charts 
\begin{align*}
\{(Exp_{\gamma_i,g},B_{2\delta_i}(0))\where \gamma_i\in H^{4,2}(S^1,S^2)
\text{ for }1\le i\le n\},
\end{align*}
and finite sets
\begin{align*}
\{W_i \in T_{\gamma_i}H^{4,2}(S^1,S^2)\where \|W_i\|_{T_{\gamma_i}H^{4,2}(S^1,S^2)}<\eps 
\text{ for }1\le i\le n\}, 
\end{align*}
there holds: If
$\alpha \in \capl_{i=1}^n Exp_{\gamma_i,g}(B_{\delta_i}(0))\subset H^{2,2}(S^1,S^2)$
satisfies
\begin{align*}
X(\alpha) = \sum_{i=1}^n \text{Proj}_{\langle W_g(\alpha)\rangle^\perp}\circ
D Exp_{\gamma_i,g}\eval_{Exp_{\gamma_i,g}^{-1}(\alpha)}(W_i)  
\end{align*}
then $\alpha$ is in $H^{4,2}(S^1,S^2)$.

\end{itemize}
\end{definition}
Property $(4)$ does not depend on the particular element $\gamma$ 
of the critical orbit $S^{1}*\gamma$, because from $\theta*X(\gamma)=X(\theta*\gamma)$ we get
\begin{align}
\label{eq:DX_theta}
D_{g}X\eval_{\gamma}=(-\theta*)\circ D_{g}X\eval_{\theta*\gamma}\circ(\theta*). 
\end{align}
and Rothe maps are invariant under conjugacy.
Concerning the regularity property $(5)$, taking $W_i=0$, we deduce that
if $X(\gamma)=0$ then $\gamma \in H^{4,2}(S^1,S^2)$. 
Furthermore, if $\gamma \in H^{4,2}(S^1,S^2)$ then the map
$\theta \mapsto \theta*\gamma$ is $C^2$ from $S^1$ to $H^{2,2}(S^{1},S^{2})$. Hence, if
$X(\gamma)=0$ then
\begin{align}
\label{eq:dot_gamma_kernel}
0 &= D_\theta(X(\theta*\gamma))\eval_{\theta=0} = D_{g}X\eval_\gamma (\dot\gamma),
\end{align}
such that the kernel of $D_{g}X\eval_\gamma$ at a critical orbit $S^1*\gamma$ is nontrivial.
The parameter $\eps>0$ ensures that $(5)$ remains stable under small perturbations
used in the Sard-Smale lemma below.
If $X$ is a vector field orthogonal to $W_{g}$ 
and $X(\gamma)=0$, then
\begin{align*}
0 &= D\big(\langle X(\alpha),W_g(\alpha)\rangle_{T_\alpha H^{2,2}(S^{1},S^{2}),g}\big)\eval_\gamma
= \langle D_{g}X\eval_\gamma,W_g(\gamma)\rangle_{T_\gamma H^{2,2}(S^{1},S^{2}),g}
\end{align*}
where the various curvature terms and terms containing derivatives of $W_g$ vanish as $X(\gamma)=0$.  
Thus, $X(\gamma)=0$ implies
\begin{align}
\label{eq:dx_to_wg_perp}
D_{g}X\eval_\gamma:\:T_\gamma H^{2,2}(S^{1},S^{2}) \to \langle W_g(\gamma)\rangle^{\perp},
\end{align}
and the projection $\text{Proj}_{\langle W_g(\gamma)\rangle^{\perp}}$ in $(4)$ is unnecessary.
\begin{lemma}
\label{l:x_kg_admissible}
The vector field $X_{k,g}$ defined in (\ref{eq:def_vector_field})
is $S^{1}$-equivariant, orthogonal to $W_g$, elliptic, and a $C^{2}$-Rothe field with respect
to the set $H^{2,2}_{reg}(S^{1},S^{2})$ of regular curves.
\end{lemma}
\begin{proof}
From Section \ref{sec:introduction}
and Section \ref{s:preliminiaries} the vector field $X_{k,g}$ is
$S^{1}$-equivariant and a $C^{2}$-Rothe field. Furthermore,
we obtain for $\alpha \in H^{2,2}(S^{1},S^{2})$
\begin{align*}
\langle X_{k,g}(\alpha)&,W_{g}(\alpha)\rangle_{T_\alpha H^{2,2}(S^{1},S^{2}),g}
= \int_{S^1}\langle \dot \alpha(t),(-D_{t,g}^{2}+1)X_{k,g}(\alpha)(t)\rangle_g\, dt\\
&= \int_{S^1}\langle \dot\alpha(t)
,-D_t\dot\alpha(t)+|\dot\alpha(t)|_{g}k(\alpha(t))J_{g}(\alpha(t))\dot\alpha(t)\rangle_g\, dt\\
&= -\int_{S^1}\langle \dot \alpha(t),D_t\dot\alpha(t)\rangle_g\, dt 
= -\int_{S^1}\frac12 \frac{d}{dt}\langle \dot \alpha,\dot\alpha\rangle_g\, dt
=0.
\end{align*} 
To show that $X_{k,g}$ is elliptic, we fix 
\begin{align*}
\{(\gamma_i,W_i) \in TH^{4,2}(S^1,S^2)\where W_i \in B_{\delta_i}(0),\, 1\le i\le n\},
\end{align*}
where $(Exp_{\gamma_i,g},B_{2\delta_i}(0))$ is a chart around $\gamma_i$, 
and $\alpha \in \capl_{i=1}^n Exp_{\gamma_i,g}(B_{\delta_i}(0))$
satisfying
\begin{align*}
X_{k,g}(\alpha) &= \sum_{i=1}^n \text{Proj}_{\langle W_g(\alpha)\rangle^\perp}\circ 
D Exp_{\gamma_i,g}\eval_{Exp_{\gamma_i,g}^{-1}(\alpha)}(W_i).  
\end{align*}
Then 
\begin{align*}
D_{t,g} \dot\alpha& - |\dot \alpha|_{g} k(\alpha)J_{g}(\alpha)\dot\alpha\\
&= (-D_{t,g}^2+1) 
\sum_{i=1}^n \text{Proj}_{\langle W_g(\alpha)\rangle^\perp}\circ 
D Exp_{\gamma_i,g}\eval_{Exp_{\gamma_i,g}^{-1}(\alpha)}(W_i).
\end{align*}
We fix $1\le i\le n$ and get
\begin{align*}
D_{t,g}^2  
&\text{Proj}_{\langle W_g(\alpha)\rangle^\perp}\circ 
D Exp_{\gamma_i,g}\eval_{Exp_{\gamma_i,g}^{-1}(\alpha)}(W_i)\\
&= 
D_{t,g}^2\big(
D Exp_{\gamma_i,g}\eval_{Exp_{\gamma_i,g}^{-1}(\alpha)}(W_i)\big)\\
&-\langle D Exp_{\gamma_i,g}\eval_{Exp_{\gamma_i,g}^{-1}(\alpha)}(W_i),W_g(\alpha)\rangle 
D_{t,g}^2W_g(\alpha),
\end{align*}
as well as
\begin{align*}
&D_{t,g}^2 \big(
D Exp_{\gamma_i,g}\eval_{Exp_{\gamma_i,g}^{-1}(\alpha)}(W_i)\big)(t)\\
&\quad= 
D^2 Exp_{\gamma_i(t),g}\eval_{Exp_{\gamma_i,g}^{-1}(\alpha)(t)} D_{t,g}^2 Exp_{\gamma_i,g}^{-1}(\alpha)(t)
(W_i(t)) 
+R_{1,i}(t)\\
&\quad= D^2 Exp_{\gamma_i(t),g}\eval_{Exp_{\gamma_i,g}^{-1}(\alpha)(t)}  
D(Exp_{\gamma_i(t),g})^{-1}\eval_{\alpha(t)}
D_{t,g}\dot\alpha (t)
(W_i(t))\\
&\qquad+R_{2,i}(t),
\end{align*}
where $R_{1,i}$ and $R_{2,i}$ consist of lower order terms containing only
derivatives of $\alpha$ up to order $1$ and
derivatives of $\gamma_i$ and $W_i$ up to order $2$. Thus
$\alpha$ is a solution of
\begin{align}
\label{eq:12}
(1-A(t))&D_{t,g}\dot\alpha = |\dot \alpha|_{g} k(\alpha)J_{g}(\alpha)\dot\alpha  +R(t) \notag\\
&
- \sum_{i=1}^n \langle D Exp_{\gamma_i,g}\eval_{Exp_{\gamma_i,g}^{-1}(\alpha)}(W_i),W_g(\alpha)\rangle 
(-D_{t,g}^2+1)W_g(\alpha), 
\end{align}
where $R$ contains only derivatives of $\alpha$ up to order $1$ and
derivatives of $\gamma_i$ and $W_i$ up to order $2$ and
$A(t)\in \mathcal{L}(T_{\alpha(t)}S^2)$ is given by
\begin{align*}
V\mapsto \sum_{i=1}^n D^2 Exp_{\gamma_i(t),g}\eval_{Exp_{\gamma_i,g}^{-1}(\alpha)(t)}  
(D(Exp_{\gamma_i(t),g})^{-1}\eval_{\alpha(t)}V)
(W_i(t)).  
\end{align*}
Since $H^{2,2}$-bounds yield $L^\infty$-bounds, choosing $\max \|W_i\|$ small enough
independently of $\{\gamma_i\}$ and $\alpha$ we may assume $\|A(t)\|<\frac12$
and $A$ is of class $H^{2,2}$ with respect to $t$. 
Since $\gamma_i$ and $W_i$ are in $H^{4,2}$ and
$(-D_{t,g}^2+1)W_g(\alpha) =\dot\alpha$,  
the right hand side of \eqref{eq:12} is in $H^{1,2}$. By standard regularity results $\alpha$ is in $H^{3,2}$,
such that the right hand side of \eqref{eq:12} is in $H^{2,2}$, which yields $\alpha \in H^{4,2}$.
Consequently, $X_{k,g}$ is elliptic.
\end{proof}

\begin{definition}
\label{d:critical_orbit}
Let $M$ be an open 
$S^1$-invariant subset of prime curves in $H^{2,2}(S^{1},S^2)$, $S^{1}*\gamma \subset M$,
and $X$ a $(M,g,S^{1})$-admissible vector field on $M$.\\  
The orbit $S^{1}*\gamma$ is called a {\em critical
orbit} of $X$, if $X(\gamma)=0$.\\ 
The orbit $S^{1}*\gamma$ is called a {\em nondegenerate critical
orbit} of $X$, if $X(\gamma)=0$ and 
\begin{align}
\label{eq:def_gamma_nondeg}
D_{g}X\eval_\gamma:\: \langle W_{g}(\gamma)\rangle^{\perp} \xpfeil{}{} \langle W_{g}(\gamma)\rangle^{\perp}   
\end{align}
is an isomorphism.\\
If $S^1*\gamma$ is critical, then using the chart $\psi_{\gamma,g}$ given in \eqref{eq:def_psi}
we define after possibly shrinking $\delta>0$
a map $X^{\gamma} \in C^2(B_\delta(0)\cap \langle W_g(\gamma)\rangle^\perp,\langle W_g(\gamma)\rangle^\perp)$
by
\begin{align}
\label{eq:def_x_gamma}
X^{\gamma}(V) := \text{Proj}_2\circ \psi_{\gamma,g}^{-1}\big(Exp_{\gamma,g}(V),X(Exp_{\gamma,g}(V))\big),
\end{align}
where $\text{Proj}_2$ denotes the projection on the second component.\\
The orbit $S^{1}*\gamma$ is called an {\em isolated critical
orbit} of $X$, if $X(\gamma)=0$ and $V=0$ is an isolated zero of $X^\gamma$.
\end{definition}
The nondegeneracy of a critical orbit does not depend on the choice of $\gamma$ in $S^{1}*\gamma$.

\begin{lemma}
\label{l:nondegenerate}
Under the assumptions of Definition \ref{d:critical_orbit}
a tangent vector $V\in B_\delta(0)\cap \langle W_g(\gamma)\rangle^\perp$ is a (nondegenerate) zero of $X^{\gamma}$ 
if and only if
$S^{1}*Exp_{\gamma,g}(V)$ is a (nondegenerate) critical orbit of $X$.
\end{lemma}
\begin{proof}
From the fact that $X(Exp_{\gamma,g}(V))\perp  W_{g}(Exp_{\gamma,g}(V))$ we get
\begin{align*}
X^{\gamma}(V)=0 \aequi X(Exp_{\gamma,g}(V))=0.  
\end{align*}
Moreover, if $X^{\gamma}(V)=0$, then
\begin{align*}
DX^{\gamma}\eval_V  &= 
\text{Proj}_2\circ D\psi_{\gamma,g}^{-1}\eval_{(Exp_{\gamma,g}(V),0)}\\
&\qquad \circ \big(DExp_{\gamma,g}\eval_{V},D_{g}X\eval_{Exp_{\gamma,g}(V)}\circ DExp_{\gamma,g}\eval_{V}\big)\\
&= A^{-1}\circ D_{g}X\eval_{Exp_{\gamma,g}(V)}\circ DExp_{\gamma,g}\eval_{V},
\end{align*}
where $A:\langle W_{g}(\gamma)\rangle^\perp\to \langle W_{g}(Exp_{\gamma,g}(V)\rangle^\perp$ is given by
\begin{align*}
A:=\text{Proj}_{\langle W_{g}(Exp_{\gamma,g}(V)\rangle^\perp}\circ D Exp_{\gamma,g}\eval_{V}.
\end{align*}
By \eqref{eq:def_trans_3} the map $A$ is an isomorphism.
Consequently, the map $DX^{\gamma}\eval_V$ 
is invertible, if and only if 
\begin{align}
\label{eq:9}
D_{g}X\eval_{Exp_{\gamma,g}(V)}\circ DExp_{\gamma,g}\eval_{V}:\: 
\langle W_{g}(Exp_{\gamma,g}(V)\rangle^\perp \xpfeil{\iso}{} \langle W_{g}(Exp_{\gamma,g}(V))\rangle^{\perp}
\end{align}
is an isomorphism.
The injectivity in \eqref{eq:9}, \eqref{eq:def_trans_1}, and (\ref{eq:dot_gamma_kernel}) implies 
that the kernel of the map
$D_{g}X\eval_{Exp_{\gamma,g}(V)}$ is given by $\langle D_tExp_{\gamma,g}(V)\rangle$. 
$D_{g}X\eval_{Exp_{\gamma,g}(V)}$ is a Rothe map and thus a Fredholm operator of index $0$, and we deduce that
\eqref{eq:9} implies the nondegeneracy of $Exp_{\gamma,g}(V)$. 
If \eqref{eq:def_gamma_nondeg} holds with $\gamma$ replaced by $Exp_{\gamma,g}(V)$, then
the kernel of $D_{g}X\eval_{Exp_{\gamma,g}(V)}$ is one dimensional, and from \eqref{eq:def_trans_1} we infer
that \eqref{eq:9} holds,
which finishes the proof.
\end{proof}


\begin{definition}
\label{d:M_g_t_s_1_homotopy}  
Let $g_t$ for $t \in [0,1]$ be a family of smooth metrics on $S^2$, which induces
a corresponding family of metrics on $H^{2,2}(S^1,S^2)$, still denoted by $g_t$. 
Let $M$ be an open 
$S^1$-invariant subset of prime curves in $H^{2,2}(S^{1},S^2)$ and $X_0$, $X_1$ 
two vector-fields on $M$ such that
$X_i$ is $(M,g_i,S^1)$-admissible for $i=0,1$.
A $C^2$ family of vector-fields $X(t,\cdot)$ on $M$ for $t\in [0,1]$ is called 
a $(M,g_t,S^1)$-homotopy between $X_0$ and $X_1$, if
\begin{itemize}
\item $X(0,\cdot)=X_0$ and $X(1,\cdot)=X_1$,
\item $\{(t,\gamma)\in[0,1]\times M \where X(t,\gamma)=0\}$ is compact,
\item $X_t:= X(t,\cdot)$ is $(M,g_t,S^1)$-admissible for all $t\in [0,1]$.
\end{itemize}
We write $(M,g,S^{1})$-homotopy, 
if the family of metrics $g_t$ is constant.
\end{definition}
If $X$ is a $(M,g_t,S^1)$-homotopy, then differentiating 
\begin{align*}
\langle X(t,\gamma),W_{g_t}(\gamma)\rangle_{T_\gamma H^{2,2}(S^1,S^2),g_t} \equiv 0  
\end{align*}
we see as
in \eqref{eq:dx_to_wg_perp} for $(t_0,\gamma_0)\in X^{-1}(0)$
\begin{align}
\label{eq:dx_to_wg_perp_homo}
D_{g}X\eval_{t_0,\gamma_0}:\rz \times T_{\gamma_0} H^{2,2}(S^1,S^2) \to \langle W_{g_{t_0}}(\gamma_0) \rangle^{\perp,g_{t_0}}. 
\end{align}
Moreover, analogous to \eqref{eq:def_x_gamma} there is $\delta>0$ such that 
\begin{align*}
X^{t_0,\gamma_0} \in C^2(B_{\delta}(t_0)\times B_{\delta}(0)\cap \langle W_{g_{t_0}}(\gamma_0) \rangle^{\perp,g_{t_0}},
\langle W_{g_{t_0}}(\gamma_0) \rangle^{\perp,g_{t_0}}),  
\end{align*}
where
\begin{align*}
X^{t_0,\gamma_0}(t,V) := \text{Proj}_3\circ \psi_{\gamma_0,t_0}^{-1}
\big(t,Exp_{\gamma_0,g_{t_0}}(V),X(t,Exp_{\gamma_0,g_{t_0}}(V))\big),  
\end{align*}
and $\psi_{\gamma_0,t_0}$ is a chart around $(t_0,\gamma_0,0)$ of the bundle
\begin{align*}
S_{[0,1]}H^{2,2}(S^{1},S^{2}) := \{&(t,\gamma,V)\in [0,1]\times TH^{2,2}(S^{1},S^{2})
\where \gamma \neq {\rm const} \\
&\text{ and } V\in \langle W_{g_t}(\gamma)\rangle^{\perp,g_t}\}, 
\end{align*}
defined in a neighborhood of $(t_0,0,0)$ in
\begin{align*}
[0,1] \times T_{\gamma_0}H^{2,2}(S^{1},S^{2}) 
\times \langle W_{g_{t_0}}(\gamma_0)\rangle^{\perp,g_{t_0}}  
\end{align*}
by
\begin{align}
\label{eq:def_psi_t0}
\psi_{\gamma,t_0}(t,V,U) &:=
\big(t,Exp_{\gamma,g_{t_0}}V,
Proj_{\langle W_{g_t}(Exp_{\gamma,g_{t_0}}V)\rangle^{\perp,g_t}}
(D Exp_{\gamma,g_{t_0}}\eval_V U)\big).
\end{align}

\begin{definition}
\label{d:nondeg_homotopy}
Let $X$ be a $(M,g_t,S^{1})$-homotopy and $(t_0,S^1*\gamma_0) \in [0,1]\times M$.  
The orbit $(t_0,S^{1}*\gamma_0)$ is called a {\em nondegenerate} zero of $X$, 
if $X(t_0,\gamma_0)=0$ and 
\begin{align}
\label{eq:def_gamma_nondeg_homotopy}
D_{g}X\eval_{(t_0,\gamma_0)}:\: 
\rz \times \langle W_{g_{t_0}}(\gamma_0)\rangle^{\perp,g_{t_0}} 
\to \langle W_{g_{t_0}}(\gamma_0)\rangle^{\perp,g_{t_0}}   
\end{align}
is surjective.
\end{definition}

Analogously to Lemma \ref{l:nondegenerate} we obtain for 
a homotopy $X$.
\begin{lemma}
\label{l:nondegenerate_homotopy}
Under the assumptions of Definition \ref{d:nondeg_homotopy}
the tuple $(t,V)$ in $B_{\delta}(t_0) \times B_\delta(0)\cap \langle W_{g_{t_0}}(\gamma)\rangle^{\perp,g_{t_0}}$ 
is a (nondegenerate) zero of $X^{t_0,\gamma_0}$ if and only if the orbit
$(t,S^{1}*Exp_{\gamma_0,g_{t_0}}(V))$ is a (nondegenerate) zero of $X$.
\end{lemma}

We give a $S^{1}$ equivariant version of the Sard-Smale lemma \cite{MR0185604,MR0322899}. 
\begin{lemma}
\label{l:s1_sard}
Let $M$ be an open 
$S^1$-invariant subset of prime curves in $H^{2,2}(S^{1},S^2)$ and $X$
a $(M,g,S^{1})$-admissible vector field on $M$. Let $\mathcal{U}$ be an open neighborhood of the zeros of $X$. 
Then there exists a $(M,g,S^{1})$-admissible vector field $Y$ such that $Y$ has only finitely many isolated,
nondegenerate zeros, $Y$ equals $X$ outside $\mathcal{U}$ and there is a $(M,g,S^{1})$-homotopy 
connecting $X$ and $Y$. 
\end{lemma}
\begin{proof}
As $X$ is proper and $X^{-1}(0)\subset H^{4,2}(S^{1},S^{2})$ 
using Lemma \ref{l:slice_lemma} we may cover $X^{-1}(0)$ with finitely many open sets
\begin{align*}
X^{-1}(0) \subset \cupl_{i=1}^{n} 
S^{1}*Exp_{\gamma_i,g}\big(B_{\delta_i}(0)\cap \langle W_g(\gamma_i)\rangle^{\perp}\big),\\
\cupl_{i=1}^{n} 
S^{1}*Exp_{\gamma_i,g}\big(B_{3\delta_i}(0)\cap \langle W_g(\gamma_i)\rangle^{\perp}\big)
\subset \mathcal{U},
\end{align*}
where $\delta_i>0$, the slice $\Sigma_{\gamma_i,g}$ is defined in $S^1\times B_{3\delta_i}(0)$, 
and $X^{\gamma_i}$ is defined in $B_{3\delta_i}(0)\cap \langle W_g(\gamma_i)\rangle^{\perp}$
for $i=1, \dots, n$.\\
Then $DX^{\gamma_i}\eval_{0}$ is in $\mathscr{R}(\langle W_g(\gamma_i)\rangle^{\perp})$, which
is open in 
$\mathscr{L}(\langle W_g(\gamma_i)\rangle^{\perp})$. 
Thus $D X^{\gamma_i}\eval_V$ remains a Rothe map for $V$ close to $0$ and
consequently a Fredholm operator of index $0$.
As Fredholm maps are locally proper,
we may assume for all $1\le i \le n$ that
the map $X^{\gamma_i}$ is proper and Rothe on
$\overline{B_{2\delta_i}(0)}\cap \langle W_g(\gamma_i)\rangle^{\perp}$, i.e.
\begin{align*}
DX^{\gamma_i}\eval_V \in \mathscr{R}(\langle W_g(\gamma_i)\rangle^{\perp}) \; \forall V \in
\overline{B_{2\delta_i}(0)}\cap \langle W_g(\gamma_i)\rangle^{\perp},\\
\overline{B_{2\delta_i}(0)}\cap \langle W_g(\gamma_i)\rangle^{\perp} \cap (X^{\gamma_i})^{-1}(K) \text{ is compact}  
\end{align*}
for all compact sets $K \subset \langle W_g(\gamma_i)\rangle^{\perp}$.\\
To construct $Y$ we proceed step by step and construct $Y_{j}$ such that
\begin{itemize}
\item[(i)] $Y_j$ equals $X$ outside 
$\cup_{i=1}^{j-1} S^{1}*Exp_{\gamma_i,g}(B_{2\delta_i}(0) \cap \langle W_g(\gamma_i)\rangle^{\perp})$,
\item[(ii)] $Y_j^{-1}(0)$ is a subset of 
\begin{align*}
\cupl_{i=1}^{n} S^{1}*Exp_{\gamma_i,g}(B_{\delta_i}(0) \cap \langle W_g(\gamma_i)\rangle^{\perp}),
\end{align*}
\item[(iii)] the critical orbits of $Y_j$ in 
  \begin{align*}
  \cupl_{i=1}^{j} S^{1}*Exp_{\gamma_i,g}(\overline{B_{\delta_i}(0)}\cap \langle W_g(\gamma_i)\rangle^{\perp}).
  \end{align*}
are isolated and nondegenerate.
\end{itemize}
Since each $X^{\gamma_i}$ is proper, $\|X(\cdot)\|$ is bounded below by a positive
constant in 
$$
\cup_{i=1}^{n} S^{1}*Exp_{\gamma_i,g}(B_{2\delta_i}(0)\setminus
\cupl_{i=1}^{n} S^{1}*Exp_{\gamma_i,g}(B_{\delta_i}(0).
$$
Consequently, $(ii)$ remains valid for all small perturbations of $X$.\\
We start with $Y_{0}:= X$.
In the $j$th step we consider $Y_{j-1}^{\gamma_j}$.
By the Sard-Smale lemma 
there is 
$V_j \in \langle W_g(\gamma_j)\rangle^{\perp}\cap T_{\gamma_j} H^{4,2}(S^{1},S^{2})$ 
arbitrarily close to zero, such
that $Y_{j-1}^{\gamma_j}+V_j$ has only nondegenerate zeros in 
$\overline{B_{\delta_i}(0)}\cap \langle W_g(\gamma_j)\rangle^{\perp}$.\\
Since $\gamma_j \in H^{4,2}(S^{1},S^{2})$, the map $\theta \mapsto \theta*\gamma_j$ is in 
$C^{2}(S^{1},H^{2,2}(S^{1},S^{2}))$ and
$S^{1}*\gamma_j$ is a $C^{2}$ sub-manifold of $H^{2,2}(S^{1},S^{2})$.
Shrinking $\delta_j>0$  we may assume the distance function $d_g(\cdot,S^{1}*\gamma_j)$
in the Riemannian manifold $H^{2,2}(S^{1},S^{2})$ satisfies
\begin{align*}
d_g(\cdot,S^{1}*\gamma_j)^2 \in C^{2}
(S^{1}*Exp_{\gamma_j,g}(B_{2\delta_j}(0)\cap \langle W_g(\gamma_j)\rangle^{\perp}),\rz),  
\end{align*}
and there are $\eps_{j,1},\eps_{j,2}>0$ such that the set  
\begin{align*}
\{\gamma \in S^{1}*Exp_{\gamma_j,g}(B_{2\delta_j}(0)\cap \langle W_g(\gamma_j)\rangle^{\perp}) 
\where \eps_{j,1}\le d_g(\gamma,S^{1}*\gamma_j)\le \eps_{j,2}\}
\end{align*}
is contained in
\begin{align*}
S^{1}*Exp_{\gamma_j,g}\big((B_{2\delta_j}(0)\setminus 
\overline{B_{\delta_j}(0)})\cap \langle W_g(\gamma_j)\rangle^{\perp}\big).
\end{align*}
We take a cut-off function $\eta \in C_c^\infty(\rz,[0,1])$ such that $\eta\equiv 1$ in 
$[0,\eps_{j,1}]$ and $\eta(x)=0$ for $x\ge\eps_{j,2}$.
Using Lemma \ref{l:slice_lemma} we define 
$$\theta_j \in C^2(S^1*Exp_{\gamma_j,g}(B_{2\delta_j}(0)\cap \langle W_g(\gamma_j)\rangle^{\perp}),S^1)$$ 
by $\theta_j:= \text{Proj}_{S^1}\circ (\Sigma_{\gamma_j,g})^{-1}$
and the vector field $Y_{j}$ on $M$ by
\begin{align*}
Y_{j}(\gamma) := Y_{j-1}(\gamma),
\end{align*}
if $\gamma \not\in S^1*Exp_{\gamma_j,g}(B_{2\delta_j}(0)\cap \langle W_g(\gamma_j)\rangle^{\perp})$ and
\begin{align*}
Y_{j}(\gamma) :=
Y_{j-1}(\gamma)+ &\eta\big(d_g(\gamma,S^{1}*\gamma_j)\big)\\
&\text{Proj}_2\circ \psi_{\theta_j(\gamma)*\gamma_j,g}
\big(Exp_{\theta_j(\gamma)*\gamma_j,g}^{-1}(\gamma),\theta_j(\gamma)*V_j\big),
\end{align*}
if $\gamma \in S^1 *Exp_{\gamma_j,g}(B_{2\delta_j}(0)\cap \langle W_g(\gamma_j)\rangle^{\perp})$.\\
Note that the map $\theta \mapsto (\theta*\gamma_j,\theta*V_j)$
is in $C^2(S^1,TH^{2,2}(S^1,S^2))$ as $(\gamma_j, V_j) \in TH^{4,2}(S^1,S^2)$. 
It is easy to see that $Y_j$ is a $S^1$ equivariant $C^2$ vector field, 
which is orthogonal to $W_g$ by construction. If $\|V_j\|$ is small enough, 
then $(i)$-$(iii)$ continue to hold for $Y_j$ as well as the Rothe property, because
Rothe maps and nondegenerate critical orbits are stable
under small perturbations.  
Moreover, $\cos(t)^2 Y_{j-1}+\sin(t)^2Y_j$ is proper for any $t\in [0,\pi/2]$, because 
$\cos(t)^2 Y_{j-1}+\sin(t)^2Y_j$ equals $Y_{j-1}$ outside 
$S^1*Exp_{\gamma_j,g}(B_{2\delta_j}(0)\cap \langle W_g(\gamma_j)\rangle^{\perp})$, which is proper,
and the zeros of $\cos(t)^2 Y_{j-1}+\sin(t)^2Y_j$ inside 
$S^1*Exp_{\gamma_j,g}(\overline{B_{2\delta_j}(0)}\cap \langle W_g(\gamma_j)\rangle^{\perp})$ 
are contained in the compact set $S^1*Exp_{\gamma_j,g}((Y_{j-1}^{\gamma_j})^{-1}([0,1]V_j))$.
If $Y_{j-1}$ is elliptic with constant $\eps_{j-1}>0$, then taking $\|V_j\|$ small enough
$\cos(t)^2 Y_{j-1}+\sin(t)^2Y_j$ remains elliptic with constant $\eps_{j}=\eps_{j-1}/2$, because
$Y_j(\gamma)$ and $Y_{j-1}(\gamma)$ differ only by
\begin{align*}
\lambda \text{Proj}_{\langle W_g(\gamma)\rangle^\perp}\circ
DExp_{\theta_j(\gamma)*\gamma_j,g}\eval_{Exp_{\theta_j(\gamma)*\gamma_j,g}^{-1}(\gamma)}
(\theta_j(\gamma)*V_j),
\end{align*}
where $\lambda \in [0,1]$ and $\theta_j(\gamma)*\gamma_j$ and $\theta_j(\gamma)*V_j$ are in $H^{4,2}$.\\
For $j=n$
we arrive at the desired vector-field $Y$.
\end{proof}
Essentially the same arguments lead to the following lemma.
\begin{lemma}
\label{l:s1_sard_ext}
Let $M$ be an open 
$S^1$-invariant subset of prime curves in $H^{2,2}(S^{1},S^2)$, $g_t$ 
for $t\in[0,1]$ a smooth family of metrics on $S^2$, and
$X$ a $(M,g_t,S^1)$-homotopy between two vector-fields  
$X_0$ and $X_1$ on $M$, which have only finitely many critical orbits in $M$,
that are all nondegenerate. 
Let $\mathcal{U}$ be an open neighborhood of the zeros of $X$. 
Then there exists a $(M,g_t,S^1)$-homotopy $Y$ and $\eps>0$ such that 
$Y_t(\gamma)=X_t(\gamma)$ for 
\begin{align*}
(t,\gamma) \in ([0,\eps] \cup [1-\eps,1])\times M \cup ([0,1]\times M)\setminus \mathcal{U},  
\end{align*}
and
\begin{align*}
DY\eval_{t,\gamma}: \rz \times \langle W_{g_t}(\gamma)\rangle^{\perp,g_t} 
\to \langle W_{g_t}(\gamma)\rangle^{\perp,g_t}   
\end{align*}
is surjective for all zeros $(t,\gamma)$ of $Y$.  
\end{lemma}
For the rest of this section we let $M$ be an open 
$S^1$-invariant subset of prime curves in $H^{2,2}(S^{1},S^2)$ and $X$
a $(M,g,S^{1})$-admissible vector field on $M$.
We shall define the $S^{1}$-equivariant Poincar\'{e}-Hopf index $\chi_{S^{1}}(X,M)$ of the vector-field $X$ with
respect to the set $M$. We begin with the definition 
of the local degree of an isolated, nondegenerate critical orbit of $X$.\\
We fix a nondegenerate critical orbit $S^1*\gamma_0$ of $X$ in $M$.
As $X$ is $(M,g,S^{1})$-admissible, $D X\eval_{\gamma_0}\in \mathcal{GR}(\langle W_g(\gamma_0)\rangle^{\perp})$
and we define the local degree of $X$ at $S^1*\gamma_0$ by
\begin{align*}
\deg_{loc,S^1}(X,S^1*\gamma_0):=\sgn D_g X\eval_{\gamma_0}.
\end{align*}
From (\ref{eq:DX_theta}) the local degree 
does not depend on the choice of $\gamma_0$ in $S^1*\gamma_0$.
\begin{definition}[$S^1$-degree]
\label{def:s1_degree}
Let $X$ be $(M,g,S^{1})$-admissible. From Lemma \ref{l:s1_sard} there is a vector field $Y$,
which is $(M,g,S^{1})$-homotopic to $X$,
with only finitely many critical orbits, that are all nondegenerate. 
The $S^1$-equivariant Poincar\'{e}-Hopf index (or $S^1$-degree) of $X$ in $M$  is defined by
\begin{align*}
\chi_{S^1}(X,M) := \sum_{\{S^{1}*\gamma\subset M \where Y(S^{1}*\gamma)=0\}}
\deg_{loc,S^1}(Y,S^{1}*\gamma).
\end{align*}
If $S^1*\gamma_0$ is an isolated critical orbit of $X$, we define the local $S^1$-degree
of $X$ in $S^1*\gamma_0$ by
\begin{align*}
\deg_{loc,S^1}(X,S^{1}*\gamma_0) := \chi_{S^1}(X,S^1*B_\delta(\gamma_0)),  
\end{align*}
where we choose $\delta>0$ such that $S^1*\gamma_0$ is the unique critical
orbit of $X$ in the geodesic ball $B_{\delta}(S^1*\gamma)$.
\end{definition}
To show that the definition does not depend on the particular choice of $Y$ or $\delta$, and
that the $S^1$-degree does not change under homotopies in the class
of $(M,g,S^{1})$-admissible vector-fields we prove

\begin{lemma}
\label{l:degree_homotopy}
Let $g_t$ for $t\in [0,1]$ be a continuous family of metrics on $H^{2,2}(S^1,S^2)$.
Suppose $X$ is a $(M,g_t,S^{1})$-homotopy between 
$X_0$ and $X_1$, such that the zeros of $X_0$ and $X_1$ are isolated and nondegenerate. 
Then
\begin{align*}
\chi_{S^1}(X_0,M)=\chi_{S^1}(X_1,M).  
\end{align*}
\end{lemma}
\begin{proof}
By Lemma \ref{l:s1_sard_ext} we may assume the the homotopy $X$ is nondegenerate, i.e.
$DX^{t,\gamma}$ is surjective whenever $X(t,S^1*\gamma)=0$.\\
Fix $(t_0,\gamma_0) \in X^{-1}(0)$. 
From the implicit function theorem, Lemma \ref{l:slice_lemma}, and Lemma \ref{l:nondegenerate_homotopy}
there is a regular $C^1$ curve $c=(c_t,c_\gamma)\in C^1(I,\rz\times M)$ 
with $I=(-1,1)$ for $t_0 \in (0,1)$ and $I=[0,1)$ for $t_0 \in \{0,1\}$,
such that $X(c(s))\equiv 0$, $c(0)=(t_0,\gamma_0)$,
and the map 
\begin{align*}
S^1 \times I \ni (\theta,s) \mapsto (c_t(s),\theta*c_\gamma(s))=\theta*c(s)  
\end{align*}
parametrizes the zero set $X^{-1}(0)$ locally around $(t_0,\gamma_0)$, 
where we define the action of $S^{1}$ on tuples $(t,\gamma)$ by
$\theta*(t,\gamma):=(t,\theta*\gamma)$.\\
The ellipticity of $X_t$ shows that $c_\gamma(s) \in H^{4,2}(S^{1},S^{2})$, thus
$\dot{c}_\gamma(s)$ is in $T_{c_\gamma(s)}H^{2,2}(S^{1},S^{2})$ and from (\ref{eq:w_g_gamma_h12})
we deduce that 
\begin{align*}
\dot{c}_\gamma(s) \text{ is transversal to } \langle W_{g_{c_t(s)}}(c_\gamma(s))\rangle^{\perp,g_{c_t(s)}}.  
\end{align*}
Since $0\neq c'(0) \in \rz \times \langle W_{g_{t_0}}(\gamma_0)\rangle^{\perp,g_{t_0}}$ 
we see from the construction of $c$ that we may assume
for all $s\in I$
\begin{align}
\label{eq:ds_c_transversal}
c'(s) \text{ is transversal to } \big(0,\dot{c}_\gamma(s)\big).
\end{align}
By the $S^{1}$-equivariance of $X$,
(\ref{eq:ds_c_transversal}), and the fact that $D_{g_{c_t(s)}}X\eval_{c_t(s),c_\gamma(s)}$ is a Fredholm operator
of index $1$ with image $\langle W_{g_{c_t(s)}}(c_\gamma(s))\rangle^{\perp,g_{c_t(s)}}$ of codimension $1$ we find
\begin{align}
\label{eq:kernel_DX_c_s}
\text{kernel}D_{g_{c_t(s)}}X\eval_{c_t(s),c_\gamma(s)} = \langle c'(s), (0,\dot{c}_\gamma(s))\rangle.  
\end{align}
Fix $(c_1,I_1)$ and $(c_2,I_2)$ such that $S^{1}*c_1(s_1)=S^{1}*c_2(s_2)$ for some $s_1\in I_1$ and $s_2 \in I_2$.
Then from the uniqueness part in the construction of $c_2$ we get $\theta_2 \in S^{1}$ such that
$\theta_2*c_2(s_2)=c_1(s_1)$.
From its construction
$\theta_2*c_2'(s_2)$ is contained in the kernel of 
$DX\eval_{c_1(s_1)}$ spanned by $\langle c_1'(s_1),(0,(\dot{c_1})_\gamma(s_1))\rangle$.
Since $c_1'(s_1)$ and $\theta_2*c_2'(s_2)$
are both transversal to $(0,(\dot{c_1})_\gamma(s_1))$ there is $0\neq\lambda_1\in \rz$ and $\lambda_2 \in \rz$
such that
\begin{align*}
\theta_2*c_2'(s_2) = \lambda_1 c_1'(s_1) + \lambda_2 (0,(\dot{c_1})_\gamma(s_1)).  
\end{align*}
We choose a function $\theta_2\in C^{1}(I,\rz/\zz)$ satisfying $\theta_2(s_2)=\theta_2$ and
$\theta_2'(s_2)=-\lambda_2$, define $\bar{c}_2\in C^{1}(I,M)$ by
$\bar{c}_2(s) := \theta_2(s)*c_2(s)$,
and get
\begin{align*}
\bar{c}_2'(s_2) &= \theta_2*c_2'(s_2)+(0,\theta_2*(\dot{c_2})_\gamma(s_2))\theta_2'(s_2)
= \lambda_1 c_1'(s_1). 
\end{align*}
With an additional change in the $s$ parameter we may easily arrive at
$\bar{c}_2'(s_2)=c_1'(s_1)$ in such a way that the map
$(\theta,s)\mapsto \theta*\bar{c}_2(s)$ still parametrizes
$S^{1}*c_2(I_2)$.
This gives a recipe how to obtain from two overlapping local parameterizations
$(c_1,I_1)$ and $(c_2,I_2)$ of $X^{-1}(0)$ a parametrization of the union 
$S^{1}*c_1(I_1)\cup S^{1}*c_2(I_2)$.
As in the classification of one dimensional manifolds \cite{MR0226651} we deduce that 
$X^{-1}(0)$
is a two dimensional manifold with components diffeomorphic to 
$S^{1}\times S^{1}$ or $S^{1}\times [0,1]$.\\ 
Let $P$ be a component of $X^{-1}(0)$ with boundary, i.e. of the type $S^{1}\times [0,1]$,
such that a parametrization of $P$ is given by
\begin{align*}
(\theta,s)\in S^{1}\times [0,1] \mapsto \theta*c(s),  
\end{align*}
where $c\in C^{1}([0,1],[0,1]\times M)$. 
First we change $c$ to arrive at
\begin{align}
\label{eq:c_in_w_perp}
c'(s) \in \rz \times \langle W_{g_{c_t(s)}}(c_\gamma(s))\rangle^{\perp,g_{c_t(s)}}
\subset \rz \times T_{c_\gamma(s)}H^{2,2}(S^{1},S^{2}).  
\end{align}
To this end we note that from the definition of $W_g$ we have
\begin{align*}
\rz \times T_{c_\gamma(s)}H^{2,2}(S^{1},S^{2})
=  
\rz \times \langle W_{g_{c_t(s)}}(c_\gamma(s))\rangle^{\perp,g_{c_t(s)}}
\oplus \langle (0,\dot{c}_\gamma(s))\rangle
\end{align*}
and denote by $\text{Proj}_1$ the projection onto   
$\rz \times \langle W_{g_{c_t(s)}}(c_\gamma(s))\rangle^{\perp,g_{c_t(s)}}$
with respect to this decomposition.
There holds
\begin{align*}
c'(s) = \text{Proj}_1(c'(s))+\lambda(s) (0,\dot{c}_\gamma(s)). 
\end{align*}
We take $\theta\in C^{1}([0,1],\rz)$ such that $\theta'(s)=-\lambda(s)$
and define $\bar{c}(s):=\theta(s)*c(s)$. Then
\begin{align*}
\bar{c}'(s) 
&= \Big(c_t'(s),\theta(s)*\big(c_\gamma'(s)-\lambda(s)\dot{c}_\gamma(s)\big)\Big)\\
&\quad \in \rz \times \theta(s)*\langle W_{g_{c_t(s)}}(c_\gamma(s))\rangle^{\perp,g_{c_t(s)}}
= \rz \times \langle W_{g_{c_t(s)}}(\bar{c}_\gamma(s))\rangle^{\perp,g_{c_t(s)}}.
\end{align*}
Thus, replacing $c$ with $\bar{c}$ we may assume (\ref{eq:c_in_w_perp}) holds.\\
Consider for $s\in [0,1]$ the family of operators 
\begin{align*}
F_s: \rz \times \langle W_{g_{c_t(s)}}(c_\gamma(s))\rangle^{\perp,g_{c_t(s)}}
\to  \rz \times \langle W_{g_{c_t(s)}}(c_\gamma(s))\rangle^{\perp,g_{c_t(s)}} 
\end{align*}
defined by
\begin{align*}
F_s(\tau,V) := \big(\langle c'(s),(\tau,V)\rangle_{\rz \times T_{c_\gamma(s)}H^{2,2}(S^{1},S^{2})},
D_{g_{c_t(s)}}X\eval_{c(s)}(\tau,V)\big).  
\end{align*}
Since 
\begin{align*}
\text{kernel}(D_{g_{c_t(s)}}X\eval_{c(s)}) \cap \rz \times \langle W_{g_{c_t(s)}}(c_\gamma(s))\rangle^{\perp,g_{c_t(s)}}
=\langle c'(s)\rangle,\\
D_{g}X\eval_{c(s)}(\rz \times \langle W_{g_{c_t(s)}}(c_\gamma(s))\rangle^{\perp,g_{c_t(s)}}) 
= \langle W_{g_{c_t(s)}}(c_\gamma(s))\rangle^{\perp,g_{c_t(s)}},
\end{align*}
each $F_s$ is an isomorphism. Moreover, the Rothe property of $X$ implies that each $F_s$ is a Rothe map,
because $F_s$ is obtained from $DX\eval_{c(s)}$ through a change in finite dimensions.
Consequently, $\sgn(F_s)$ is well defined and by its homotopy invariance independent of $s\in [0,1]$. 
If $c_t'(s) \neq 0$ we have again by the homotopy invariance
$\sgn(F_s)=\sgn(\tilde{F}_s)$, where
\begin{align*}
\tilde{F}_s(\tau,V):= F_s\big(\tau,V+(c_t'(s))^{-1}\tau c_\gamma'(s)\big).  
\end{align*}
We have
\begin{align*}
\tilde{F}_s &=
\begin{pmatrix}
(c_t'(s))^{-1}\|c'(s)\|^{2} & \langle c_\gamma'(s),\cdot\rangle\\
0                           & D_\gamma X\eval_{c(s)}  
\end{pmatrix}
\sim
\begin{pmatrix}
(c_t'(s))^{-1}\|c'(s)\|^{2} & 0\\
0                           & D_\gamma X\eval_{c(s)}  
\end{pmatrix}.
\end{align*}
Hence, for all $s\in [0,1]$ such that $c_t'(s) \neq 0$ there holds
\begin{align}
\label{eq:sgn_F_s_sgn_D_gamma}
\sgn(F_s) &= \sgn(\tilde{F}_s)=\sgn(c_t'(s)) \sgn(D_\gamma X\eval_{c(s)}).  
\end{align}
Let $S^1*\alpha_1,\dots,S^1*\alpha_{k_0}$ be the critical orbits of $X_0$ and
$S^1*\beta_1,\dots,S^1*\beta_{k_1}$ be the critical orbits of $X_1$.
The critical orbits of $X_0$ and $X_1$ are boundary points of $X^{-1}(0)$.
From (\ref{eq:sgn_F_s_sgn_D_gamma}) we get
\begin{itemize}
\item $\sgn D X_0\eval_{\alpha_i}=-\sgn D X_0\eval_{\alpha_j}$, if $S^1*\alpha_i$ and $S^1*\alpha_j$
are boundary orbits of the same component of $X^{-1}(0)$,
\item $\sgn D X_1\eval_{\beta_i}=-\sgn D X_1\eval_{\beta_j}$, if $S^1*\beta_i$ and $S^1*\beta_j$
are boundary orbits of the same component of $X^{-1}(0)$,
\item $\sgn D X_0\eval_{\alpha_i}=\sgn D X_1\eval_{\beta_j}$, if $S^1*\alpha_i$ and $S^1*\beta_j$
are boundary orbits of the same component of $X^{-1}(0)$.
\end{itemize}
Putting the above facts together, we see that
\begin{align*}
\chi_{S^1}(X_0,M)=\chi_{S^1}(X_1,M).  
\end{align*}
\end{proof}

\section{The Degree of an Isolated Critical Orbit}
\label{sec:deg_poincare}
Let $\gamma \in H^{2,2}(S^1,S^2)$ be a prime, regular curve
such that $S^1*\gamma$ is an isolated critical orbit of $X_{k,g}$.
Then the curve
\begin{align*}
\mu(t) := \gamma\big(t|\dot\gamma|_g^{-1}\big)  
\end{align*}
is a closed $k$-magnetic geodesic with minimal period $\omega:= |\dot\gamma|_g$
such that
$t \mapsto (\mu(t),\dot\mu(t))$ lies
in the bundle $$E_1:=\{(x,V) \in TS^2\where |V|_g=1\}.$$
We fix a transversal section $\Sigma$ 
in $E_1$ at the point $\theta:=(\gamma(0),\omega^{-1}\dot\gamma(0))$
and denote by $P:B_1\cap \Sigma \to B_2 \cap \Sigma$ 
the corresponding Poincar\'{e} map, where $B_1,\,B_2$ are open neighborhoods
of $\theta$
 (see \cite[Chap. 7-8]{MR515141}).\\
In this section we shall show 
\begin{lemma}
\label{l:degree_poincare}
Under the above assumptions $\theta$ is an isolated fixed point
of $P$ and there holds
\begin{align*}
\deg_{loc, S^{1}}(X_{k,g},S^1*\gamma)= -i(P,\theta)\ge -1,  
\end{align*}
where $i(P,\theta)$ denotes the index of the isolated fixed point $\theta$.  
\end{lemma}
We consider the linearizations of equation \eqref{eq:magnetic_geodesic}
and \eqref{eq:1} given by
\begin{align}
\label{eq:jacobi_magnetic}
0&= -D_{t,g}^2 V - R_g\big(V,\dot\mu \big)\dot\mu +k(\mu)J_{g}(\mu) D_{t,g}V \notag\\
&\quad +\big(k'(\mu)V\big)J_{g}(\mu)\dot\mu
+ k(\mu) \big(D_gJ_{g}\eval_{\mu}V\big)\dot\mu.  
\end{align}
and
\begin{align}
\label{eq:jacobi_prescribed}
0&= -D_{t,g}^2 V - R_g\big(V,\dot\gamma \big)\dot\gamma
+|\dot\gamma|_g k(\gamma)J_{g}(\gamma) D_{t,g}V
 \notag \\
&\quad +|\dot\gamma|_g \big(k'(\gamma)V\big)J_{g}(\gamma)\dot\gamma
+ |\dot\gamma|_g k(\gamma) \big(D_gJ_{g}\eval_{\gamma}V\big)\dot\gamma \notag\\
&\quad+ |\dot\gamma|_g^{-1}
\langle D_{t,g}V,\dot\gamma\rangle_g k(\gamma)J_{g}(\gamma)\dot\gamma \notag \\
&= -D_{t,g}^2 V +T(V,D_{t,g}V),  
\end{align}
where $T(V,D_{t,g}V)$ abbreviates all terms containing $V$ or $D_{t,g}V$.  
For $(V_1,V_2)$ in $T_{\gamma(0)}S^{2}\times T_{\gamma(0)}S^{2} $ we denote
by $\Phi(\cdot,(V_1,V_2))$, respectively $U(\cdot,(V_1,V_2))$, 
the solution to \eqref{eq:jacobi_magnetic}, respectively \eqref{eq:jacobi_prescribed}, 
with initial values
\begin{align*}
V(0)=V_1 \text{ and } D_{t,g}V(0)=V_2.  
\end{align*}
Then
\begin{align*}
dP\eval_\theta = \text{Proj}_{T_\theta\Sigma} \circ 
\big(\Phi(\omega,\cdot),D_{t,g}\Phi(\omega,\cdot)\big)\eval_{T_\theta\Sigma},  
\end{align*}
where $\text{Proj}_{T_\theta\Sigma}$ is 
the projection onto $T_\theta \Sigma$ with kernel given by
\begin{align*}
\langle(\dot\mu(0),D_{t,g}\dot\mu(0))^{T},(0,\dot\mu(0))^{T}\rangle.  
\end{align*}

\begin{lemma}
\label{l:degree_poincare_nondeg}
Suppose $\theta$ is a nondegenerate fixed point of $P$, i.e.
the linearized Poincar\'{e} map $dP\eval_\theta:\: T_\theta \Sigma \to T_\theta \Sigma$
has no eigenvalues equal to one. Then $S^{1}*\gamma$ is a nondegenerate 
critical orbit of $X_{k,g}$ and
\begin{align*}
\deg_{loc, S^{1}}(X_{k,g},S^1*\gamma)=-\sgn\big(\det(dP\eval_\theta-I)\big).  
\end{align*}
\end{lemma}
\begin{proof}
Since index and nondegeneracy do not depend on the transversal section we may assume
$T_\theta \Sigma = T_\omega \Sigma$, where we write for $q \in \rz$
\begin{align*}
T_q \Sigma := \{&(V_1,V_2)^{T}\in 
T_\theta(TS^{2}) \iso T_{\gamma(0)}S^{2}\times T_{\gamma(0)}S^{2} \where\\
&(V_1,V_2)^{T} \text{ is orthogonal to } 
(\dot\gamma(0),q D_{t,g}\dot\gamma(0))^{T} \text{ and }
(0,\dot\gamma(0))^{T}\}   
\end{align*}
with respect to the componentwise scalar product.\\
From (\ref{eq:1}) and the symmetries
of the curvature tensor and $J_g$ we obtain for 
any solution $V$ to \eqref{eq:jacobi_magnetic} or \eqref{eq:jacobi_prescribed} 
\begin{align}
\label{eq:energy_constant}
\frac{d}{dt}\langle D_{t,g}V(t),\dot\gamma(t)\rangle_{g}
&= \langle D_{t,g}^2V(t),\dot\gamma(t)\rangle_{g}+
\langle D_{t,g}V(t),D_{t,g}\dot\gamma(t)\rangle_{g} \notag \\
&= \langle |\dot \gamma|_g k(\gamma)(D_g J_g\eval_{\gamma(t)} V(t))\dot\gamma(t)
,\dot\gamma(t)\rangle_{g}\notag \\
&=0,
\end{align}
because for a variation $(s,t) \mapsto \Gamma(s,t)$ of $\gamma$ with 
$\rand_s \Gamma(0,\cdot)=V$ we find
\begin{align*}
0 &= \frac{d}{ds}\langle J_g(\Gamma(s,t))\rand_t\Gamma(s,t)
,\rand_t\Gamma(s,t)\rangle_g\eval_{s=0}\\
&=\langle (D_gJ_g\eval_{\gamma(t)}V(t))\dot\gamma(t),\dot\gamma(t)\rangle_g\\
&\quad +\underbrace{\langle J_g(\gamma(t))V(t),\dot\gamma(t)\rangle_g 
+\langle J_g(\gamma(t))\dot\gamma(t),V(t)\rangle_g}_{=0}    
\end{align*}
Up to scaling with $\omega$ in $t$
equations \eqref{eq:jacobi_magnetic} and \eqref{eq:jacobi_prescribed}
only differ by
\begin{align*}
|\dot\gamma|_g^{-1}\langle D_{t,g}V,\dot\gamma\rangle_g k(\gamma)J_{g}(\gamma)\dot\gamma.
\end{align*}
Since $\langle D_{t,g}U,\dot\gamma\rangle_g$ is constant we get
for $V_2$ orthogonal to $\dot\gamma(0)$
\begin{align*}
U\big(t \omega^{-1},(V_1,\omega V_2)\big)= \Phi\big(t,(V_1,V_2)\big).  
\end{align*}
Consequently
\begin{align*}
dP\eval_\theta &= \text{Proj}_{T_\omega\Sigma} \circ A_{\omega}^{-1}\circ 
\big(U(1,\cdot),D_{t,g}U(1,\cdot)\big)\circ A_{\omega}\eval_{T_\omega \Sigma}\\
\end{align*}
where $A_{\omega},\, \text{Proj}_{T_\theta\Sigma} 
\in \mathscr{L}(T_{\gamma(0)}S^{2}\times T_{\gamma(0)}S^{2})$
are given by
\begin{align*}
A_{\omega}(V_1,V_2)&:= (V_1, \omega V_2),\\
\text{Proj}_{T_\theta\Sigma}\begin{pmatrix}V_1\\V_2\end{pmatrix}&:=
\begin{pmatrix}V_1\\V_2\end{pmatrix}
-\frac
{\big\langle\begin{pmatrix}V_1\\V_2\end{pmatrix},
\begin{pmatrix}\dot\gamma(0)\\\omega D_{t,g}\dot\gamma(0)\end{pmatrix}\big\rangle_g}
{\big\|\begin{pmatrix}\dot\gamma(0) \\D_{t,g}\dot\gamma(0)\end{pmatrix}\big\|_g^{2}}
\begin{pmatrix}\dot\gamma(0)\\ \omega^{-1} D_{t,g}\dot\gamma(0)\end{pmatrix}.  
\end{align*}
We note that $A_\omega:\: T_\theta\Sigma \xpfeil{\iso} 
T_1\Sigma$ and
$
\text{Proj}_{T_\theta\Sigma} \circ A_{\omega}^{-1} =
A_{\omega}^{-1} \circ 
\text{Proj}_{T_1\Sigma}^{\perp}.
$
Hence, we may replace in the following 
$dP\eval_\theta$ by $dP\eval_{\gamma}:T_1\Sigma \to T_1\Sigma$
\begin{align*}
dP\eval_{\gamma} :=
\text{Proj}_{T_1\Sigma}^{\perp}
\circ 
\big(U(1,\cdot),D_{t,g}U(1,\cdot)\big)
\eval_{T_1\Sigma}.
\end{align*} 
To show that $S^{1}*\gamma$ is a nondegenerate 
critical orbit we
fix $V \in \langle W_g(\gamma)\rangle^\perp$
such that $DX_{k,g}\eval_\gamma (V)=0$. There are $\lambda_1,\lambda_2 \in \rz$ such that
\begin{align*}
\begin{pmatrix}W_1\\W_2\end{pmatrix}
:= \begin{pmatrix}V(0)\\ D_{t,g}V(0)\end{pmatrix}
+ \lambda_1 \begin{pmatrix}\dot\gamma(0)\\ D_{t,g}\dot\gamma(0)\end{pmatrix}
+\lambda_2 \begin{pmatrix}0\\ \dot\gamma(0)\end{pmatrix}
\in T_1\Sigma.
\end{align*}
Using the fact that $V$, $\dot\gamma$, and $t\dot\gamma$ solve (\ref{eq:jacobi_prescribed})
we get
\begin{align*}
\text{Proj}_{T_1\Sigma}^{\perp}
\begin{pmatrix}U(1,(W_1,W_2))\\D_{t,g}U(1,(W_1,W_2))\end{pmatrix}
&= \begin{pmatrix}W_1\\W_2\end{pmatrix}.
\end{align*}
Since $dP\eval_\gamma$ has no eigenvalues equal to one $(W_1,W_2)$ equals $(0,0)$
and $V=\lambda_1 \dot\gamma + \lambda_2 t \dot\gamma$. From the periodicity of $V$
we obtain $\lambda_2=0$, and the fact that  $V \in \langle W_g(\gamma)\rangle^\perp$ 
gives $\lambda_1=0$. Consequently, $S^{1}*\gamma$ is a nondegenerate critical orbit
of $X_{k,g}$.\\
We consider $\tilde{X} \in \mathscr{L}(T_\gamma H^{2,2}(S^1,S^2)\times \rz)$
defined by
\begin{align*}
\tilde{X}(V,\delta) := \big(DX_{k,g}\eval_\gamma(V)+\delta W_g(\gamma)
,\delta-\eps\langle V,W_g(\gamma)\rangle_g\big),  
\end{align*}
where $\eps>0$ will be chosen later.
Then $\tilde{X}$ is of the form $identity + compact$. 
If $\tilde{X}(V_0,\delta_0)=0$, then $\delta_0=0$, since
$DX_{k,g}\eval_\gamma(V_0)$ is orthogonal to $W_g(\gamma)$. 
From $\eps>0$ we get
$V_0\in \langle W_g(\gamma)\rangle^\perp$ 
and finally $V_0=0$. Thus, $\tilde{X}$ is invertible.
With respect to the decomposition
\begin{align*}
T_\gamma H^{2,2}(S^1,S^2)\times \rz 
= \langle W_g(\gamma)\rangle^\perp \times \{0\} 
\oplus \langle W_g(\gamma)\rangle \times \{0\}
\{0\} \times \rz    
\end{align*}
we have
\begin{align*}
\tilde{X}=
\begin{pmatrix}
DX_{k,g}\eval_\gamma & * & 0\\
0 &0 &1\\
0 &-\eps &1    
\end{pmatrix},
\end{align*}
such that
\begin{align*}
\deg(\tilde{X},(0,0)) =
\deg(DX_{k,g}\eval_\gamma\eval_{\langle W_g(\gamma)\rangle^\perp},0)
= \deg_{loc, S^{1}}(X_{k,g},S^1*\gamma),   
\end{align*}
where $\deg$ denotes the usual Leray-Schauder degree.\\
Using the ideas in \cite[Chap. 3]{MR736839} we define a homotopy 
\begin{align*}
\Phi: [0,1]\times T_\gamma H^{2,2}(S^1,S^2)\times \rz \to T_\gamma H^{2,2}(S^1,S^2)\times \rz,  
\end{align*}
by
\begin{align*}
\Phi&(s,(V,\delta)):=\\
&\bigg(V+s(-D_{t,g}^2+1)^{-1}\big(T(V,D_{t,g}V)-V\big) + s\delta W_g(\gamma)\\
&\quad +(1-s)(-D_{t,g}^2+1)^{-1}\\
&\quad \Big(D_{t,g}^2 U(\cdot,(V(0),D_{t,g}V(0)),\delta)
- U(\cdot,(V(0),D_{t,g}V(0)),\delta)\Big),\\
&\quad\delta -\eps 
\int_0^1 (-D_{t,g}^2+1)\Big((1-s)V+sU(\cdot,(V(0),D_{t,g}V(0)),\delta)\Big)\dot\gamma\bigg),  
\end{align*}
where 
$(-D_{t,g}^2+1)^{-1}$ maps $T_\gamma L^{2}(S^{1},S^{2})$ to $T_\gamma H^{2,2}(S^{1},S^{2})$
and the function $U(\cdot,(V_1,V_2),\delta)$ denotes the solution to
\begin{align*}
0= -D_{t,g}^{2}V + T(V,D_{t,g}V)+\delta \dot\gamma,  
\end{align*}
with initial values $V(0)=V_1$ and $D_{t,g}V(0)=V_2$.\\
Fix $(s_0,V_0,\delta_0) \in [0,1]\times T_\gamma H^{2,2}(S^1,S^2)\times \rz$ such that
\begin{align*}
\Phi(s_0,(V_0,\delta_0))=(0,0).  
\end{align*}
Then $V_0$ is a periodic solution of
\begin{align}
\label{eq:homotopy_soln}
0&= -D_{t,g}^{2}V_0+V_0+sT(V_0,D_{t,g}V_0)-sV_0+s\delta \dot\gamma +(1-s_0)\notag\\  
&\quad\Big(D_{t,g}^2 U(\cdot,(V(0),D_{t,g}V(0)),\delta)
- U(\cdot,(V(0),D_{t,g}V(0)),\delta)\Big).
\end{align}
Since $U(\cdot,(V_0(0),D_{t,g}V_0(0)),\delta_0)$ is a solution to (\ref{eq:homotopy_soln})
with the same initial values we see that
\begin{align*}
V_0= U(\cdot,(V_0(0),D_{t,g}V_0(0)),\delta_0).  
\end{align*}
In this case $\Phi(s_0,(V_0,\delta_0))=(0,0)$ is equivalent to
\begin{align*}
\big(DX_{k,g}\eval_\gamma(V_0)+\delta_0 W_g(\gamma)
,\delta_0-\eps\langle V_0,W_g(\gamma)\rangle_g\big) =(0,0),
\end{align*}
which shows that $V_0=0$ and $\delta_0=0$. Consequently,
\begin{align*}
\deg(\tilde{X},(0,0)) &= \deg(\Phi(1,\cdot),(0,0))
= \deg(\Phi(0,\cdot),(0,0)).
\end{align*}
We choose $\tilde{E}_i=(E_i,\delta_{i}) \in T_\gamma H^{2,2}(S^1,S^2)\times \rz$ 
for $1\le i \le 5$ such that 
\begin{align*}
\{(E_i(0),D_{t,g}E_i(0),\delta_{i})^{T}\where 1 \le i \le 5\}
\end{align*}
is an orthonormal basis
of $T_{\gamma(0)}S^{2}\times T_{\gamma(0)}S^{2} \times \rz$
with respect to the componentwise scalar product.
Since $\Phi(0,(V,\delta))=(V,\delta)$ for all 
\begin{align*}
(V,\delta) \in W_0:=
\{(V,0) \in T_\gamma H^{2,2}(S^1,S^2) \times \rz\where
V(0)=0=D_{t,g}V(0)\},   
\end{align*}
there holds for $W_1:= \langle \tilde{E}_i:\: 1 \le i \le 5 \rangle$
\begin{align*}
\deg(\Phi(0,\cdot),(0,0)) = \deg(P_{W_1}\circ \Phi(0,\cdot)\eval_{W_1},0),   
\end{align*}
where
$P_{W_1}: T_\gamma H^{2,2}(S^1,S^2)\times \rz \to W_1$ is given by
\begin{align*}
P_{W_1}(V,\delta) := 
\sumg_{i=1}^{5} 
\big(\langle V(0),E_i(0)\rangle+ \langle D_{t,g} V(0), D_{t,g} E_i(0)\rangle
+\delta\delta_i\big)  
\tilde{E}_i.
\end{align*}
Note that $P_{W_1}$ is the projection onto $W_1$ with kernel $W_0$.\\
We define $Ev_0:\: T_\gamma H^{2,2}(S^1,S^2)\times \rz
\to T_{\gamma(0)}S^2 \times T_{\gamma(0)}S^2\times \rz$ by
\begin{align*}
Ev_0(V,\delta) := \big(V(0),D_{t,g}V(0),\delta\big).  
\end{align*}
Then $Ev_0\eval_{W_1}$ is an isomorphism and we have
\begin{align*}
\deg(P_{W_1}\circ \Phi(0,\cdot)\eval_{W_1},0)
=
\deg(Ev_0 \circ P_{W_1}\circ \Phi(0,\cdot)\circ (Ev_0\eval_{W_1})^{-1},0).  
\end{align*}
We note that for a function $U$
\begin{align*}
(-D_{t,g}^2+1)^{-1}
\big((-D_{t,g}^2+1) U\big)
= U + Q,  
\end{align*}
where $Q$ solves $(-D_{t,g}^2+1)Q=0$ with boundary
conditions
\begin{align*}
Q(0)-Q(1) = U(1)-U(0),\\
D_{t,g}Q(0)-D_{t,g}Q(1) = D_{t,g}U(1)-D_{t,g}U(0).
\end{align*}
We let $B_1$ and $B_2$ be the smooth parallel vector fields along $\gamma$
such that $\{B_1(0),B_2(0)\}$ is a basis of $T_{\gamma(0)}S^2$.
Then the set $\Lambda_0$ of functions $Q$ with $(-D_{t,g}^2+1)Q=0$ is given by
\begin{align*}
\Lambda_0 := \{e^{t}\sum_{i=1}^2\lambda_i B_i(t)+ e^{-t}\sum_{i=1}^2\mu_i B_i(t)\where
\lambda_1,\,\lambda_2,\,\mu_1,\,\mu_2 \in \rz\}. 
\end{align*}
We define $L_0,\,L_1:\: \Lambda_0 \to T_{\gamma(0)}S^2 \times T_{\gamma(0)}S^2$
by
\begin{align*}
L_0(Q) := \big(Q(0)-Q(1),D_{t,g}Q(0)-D_{t,g}Q(1)\big),\\
L_1(Q) := \big(Q(0),D_{t,g}Q(0)\big).  
\end{align*}
It is easy to see that $L_0$ and $L_1$ are isomorphisms. We have
\begin{align*}
Ev_0 &\circ P_{W_1}\circ \Phi(0,\cdot)\circ (Ev_0\eval_{W_1})^{-1}(V_1,V_2,\delta)\\  
&= 
\bigg(
L_1\circ L_0^{-1}
\begin{pmatrix}
U(1,(V_1,V_2),\delta)-U(0,(V_1,V_2),\delta)\\
D_{t,g}U(1,(V_1,V_2),\delta)-D_{t,g}U(0,(V_1,V_2),\delta)  
\end{pmatrix}
,\\
&\quad \delta -\eps
\int_0^1 (-D_{t,g}^2+1)\Big(U(\cdot,(V_1,V_2),\delta)\Big)\dot\gamma\bigg).
\end{align*}
The map $L_1\circ L_0^{-1}$ may be computed explicitly 
solving a system of linear equations.
Using the fact that the parallel transport is an isometry it is easy to see  
that $\det(L_1\circ L_0^{-1})>0$. Thus we may replace $L_1\circ L_0^{-1}$
by $id$ without changing the degree.
Hence we need to compute the degree of 
$Y \in \mathscr{L}(T_{\gamma(0)}S^2 \times T_{\gamma(0)}S^2\times \rz)$
given by
\begin{align*}
Y(V_1,V_2,\delta) &:=
\bigg(
\begin{pmatrix}
U(1,(V_1,V_2),\delta)-U(0,(V_1,V_2),\delta)\\
D_{t,g}U(1,(V_1,V_2),\delta)-D_{t,g}U(0,(V_1,V_2),\delta)  
\end{pmatrix}
,\\
&\quad \delta -\eps
\int_0^1 (-D_{t,g}^2+1)\Big(U(\cdot,(V_1,V_2),\delta)\Big)\dot\gamma\bigg). 
\end{align*}
To compute $\deg(Y,0)$ we decompose
$T_{\gamma(0)}S^2 \times T_{\gamma(0)}S^2\times \rz$ into
\begin{align*}
T_1\Sigma\times \{0\} 
\oplus 
\big\langle(\dot\gamma(0),D_{t,g}\dot\gamma(0),0)\big\rangle  
\oplus 
\big\langle(0,\dot\gamma(0),0)\big\rangle  
\oplus 
\big\langle(0,0,1)\big\rangle,  
\end{align*}
where the decomposition is orthogonal with respect to the componentwise scalar product.
We have
\begin{align*}
Y(\dot\gamma(0),D_{t,g}\dot\gamma(0),0)
&= \big(0,0,-\eps \|\dot\gamma\|^2_{H^{1,2}}\big),\\
Y(0,\dot\gamma(0),0)
&= \big(\dot\gamma(0),D_{t,g}\dot\gamma(0),
-\eps (\|\sqrt{t}\dot\gamma\|^2_{L^2}
+\|\sqrt{t}D_{t,g}\dot\gamma\|^2_{L^2})\big).  
\end{align*}
We obtain analogously to \eqref{eq:energy_constant} 
\begin{align*}
\frac{d}{dt}\langle D_{t,g}U(t,(0,0),1),\dot\gamma(t)\rangle_g
= \langle \dot\gamma(t),\dot\gamma(t)\rangle_g>0,\\
\langle Y(0,0,1),(0,\dot\gamma(0),0)\rangle_g 
= \langle D_{t,g}U(1,(0,0),1),\dot\gamma(0)\rangle_g >0.
\end{align*}
Choosing $\eps>0$ small enough we find
\begin{align*}
\langle Y(0,0,1),(0,0,1)\rangle_g = 1-O(\eps)>0.
\end{align*}
Moreover, using again \eqref{eq:energy_constant} we get
for all $(V_1,V_2) \in T_1\Sigma$
\begin{align*}
\langle Y(V_1,V_2,0),(0,\dot\gamma(0),0)\rangle_g =0.
\end{align*}
Consequently we obtain with respect to the above decomposition
\begin{align*}
\deg(Y,0) &=\sgn \det
\begin{pmatrix}
P\eval_\gamma - id & 0 & 0 & *\\
* & 0 & 1 & *\\
0 & 0 & 0 & +\\
* & - & - & +  
\end{pmatrix}\\
&= -\sgn\det (P\eval_\gamma - id),
\end{align*}
which proves the claim. 
\end{proof}

\begin{proof}[Proof of Lemma \ref{l:degree_poincare}]
The fact that $\theta$ is an isolated fixed point is obvious from
the properties of the Poincar\'{e} map $P$ (see \cite[Thm. 7.1.2]{MR515141}).\\
We may choose $\delta>0$ such that $S^1*\gamma$ is the unique critical
orbit of $X_{k,g}$ in the geodesic ball $B_{\delta}(S^1*\gamma)$,
$\theta$ is the unique fixed point of $P$ in the geodesic ball $B_{\delta}(\theta)$,
$i(P,\theta)= i(P,B_{\delta}(\theta))$, and
\begin{align*}
\deg_{S^1}(X_{k,g},B_{\delta}(S^1*\gamma))= \deg_{S^1,loc}(X_{k,g},S^1*\gamma).    
\end{align*}
From the homotopy invariance and 
a Kupka-Smale theorem for magnetic flows in \cite{MR2250797}
we may assume by Lemma \ref{l:degree_poincare_nondeg} that the critical orbits of $X_{k,g}$ 
in $B_{\delta}(S^1*\gamma)$ and the fixed points of $P$ in $B_\delta(\theta)$
are nondegenerate. Using again Lemma \ref{l:degree_poincare_nondeg} we find
\begin{align*}
\deg_{S^1}(X_{k,g},B_{\delta}(S^1*\gamma)) = -i(P,B_{\delta}(\theta)).  
\end{align*}
Finally, since we may assume that 
the Poincar\'{e} map is area preserving (see \cite[Thm. 8.1.3]{MR515141}), we obtain
from \cite{MR0353372,MR0346846}
\begin{align*}
i(P,\theta)\le 1.  
\end{align*}
This yields the claim.
\end{proof}
\section{The Unperturbed Problem}
\label{sec:unperturbed}
Let $S^{2}=\rand B_1(0) \subset \rz^{3}$ be the standard round sphere with induced metric $g_0$.
Then the prescribed geodesic curvature equation with $k\equiv k_0$ on $(S^2,g_0)$ is given by
\begin{align}
\label{eq:unpert:1}
Proj_{\gamma^\perp}\ddot \gamma = |\dot \gamma| k_0 \gamma \cross \dot \gamma,  
\end{align}
where $\gamma \in H^{2,2}(S^1,S^2)$, $\dot \gamma$ and $\ddot \gamma$ are the usual derivatives
of $\gamma$ considered as a curve in $\rz^3$, $|\dot \gamma|$ 
is the Euclidean norm of $\dot \gamma$
in $\rz^3$.\\
To compute the $S^1$-degree of the unperturbed equation \eqref{eq:unpert:1}
we proceed in three steps. {\em Step 1:} We compute explicitly the set $\mathcal{Z}$
of simple solutions in $H^{2,2}(S^1,S^2)$ to \eqref{eq:unpert:1}
and show that $\mathcal{Z}$ is a finite dimensional, nondegenerate manifold, 
in the sense that we have for all $\alpha \in \mathcal{Z}$
\begin{align*}
T_\alpha \mathcal{Z} = \text{kernel}(D_{g_0}X_{k_0,g_0}\eval_\alpha),\\
T_\alpha H^{2,2}(S^1,S^2) = T_\alpha \mathcal{Z} \oplus R(D_{g_0}X_{k_0,g_0}\eval_\alpha).  
\end{align*}
{\em Step 2:} We perform a finite dimensional reduction of a slightly perturbed problem:
We consider for $k_1\in C^2(S^2,\rz)$, which will be chosen later,
and $\eps \in \rz$, which is assumed to be very small,
the perturbed vector field $X_{g_0,\eps}$ defined by
\begin{align*}
X_{g_0,\eps}(\gamma)&:=
(-D_{t,g_0}^{2} + 1)^{-1} 
\big(-D_{t,g_0} \dot\gamma + |\dot \gamma|_{g_0} 
(k_0+\eps k_1(\gamma))\gamma \times \dot\gamma\big)\\
&= X_{k_0,g_0}(\gamma)+\eps K_1(\gamma),
\end{align*}
where the vector field $K_1$ is given by
\begin{align*}
K_1(\gamma):= (-D_{t,g_0}^{2} + 1)^{-1} |\dot \gamma|_{g_0} 
\big(k_1(\gamma) \gamma \times \dot\gamma\big).  
\end{align*}
We show that if $\alpha_0 \in \mathcal{Z}$ is a nondegenerate zero
of the vector field $\alpha \mapsto P_1(\alpha) \circ K_1(\alpha)$
on $\mathcal{Z}$, 
where $P_1(\alpha)$ is a projection onto $T_\alpha \mathcal{Z}$ defined below,
then there is a unique nondegenerate critical orbit $S^1*\gamma(\eps)$ for any
$0<\eps <<1$ such that $\gamma(\eps)$ converges to $\alpha_0$ as $\eps \to 0^+$
and
\begin{align*}
\deg_{loc, S^{1}}(X_{g_0,\eps},S^1*\gamma(\eps)) 
&= -\deg_{loc}(P_1(\cdot)\circ K_1(\cdot),\alpha_0).  
\end{align*}  
{\em Step 3:} We choose
\begin{align*}
k_1(x):= \langle x,e_3\rangle \text{ for } x \in S^{2}=\rand B_1(0) \subset \rz^{3},  
\end{align*}
where $\{e_1,e_2,e_3\}$ denotes the standard basis of $\rz^{3}$
and show that $P_1(\cdot)\circ K_1(\cdot)$ has exactly two nondegenerate zeros
of degree $+1$. This yields the formula
$\chi_{S^1}(X_{k_0,g_0},M)=-2$, where $M$ is the subset of 
$H^{2,2}(S^{1},S^{2})$ consisting of simple and regular curves.
\subsection{The simple solutions of \eqref{eq:unpert:1}}
\label{subsec:soln_unperturbed}
Differentiating
twice the identity $|\gamma|^2=1$ we find 
$\langle \ddot \gamma,\gamma\rangle + |\dot \gamma|^2 \equiv 0$
and (\ref{eq:unpert:1}) is equivalent to
\begin{align}
\label{eq:unpert:2}
\ddot \gamma &= |\dot \gamma| k_0 \gamma \cross \dot \gamma - |\dot \gamma|^2 \gamma.   
\end{align}
In order to solve the ordinary differential (\ref{eq:unpert:2}) we fix initial conditions
\begin{align*}
\gamma(0)=\gamma_0 \in S^2 \text{ and } \dot\gamma(0)=\tilde{v}_0 \in T_{\gamma_0}S^2. 
\end{align*}
If $\tilde{v}_0=0$ then $\gamma$ is given by the constant curve $\gamma \equiv \gamma_0$.
We may assume in the sequel
\begin{align*}
\lambda:= |\tilde{v}_0|>0.  
\end{align*}
If $k_0\neq 0$ then there is a unique $r=r(k_0) \in (-1,1)\setminus \{0\}$ such that
\begin{align*}
k_0= \frac{\sqrt{1-r^2}}{r}.  
\end{align*}
For $k_0=0$, the case of geodesics, we may take $r=\pm 1$.\\
For $\lambda>0$ and a positive oriented orthonormal system $\{v_0,v_1,w\}$ we define
the function $\alpha \in C^\infty(\rz,S^2)$ by
\begin{align*}
\alpha(t,\lambda,v_0,v_1,w) := \sqrt{1-r^2}w+r\cos(\lambda r^{-1}t)v_1+r\sin(\lambda r^{-1}t)v_0  
\end{align*}
A direct calculation shows that $\alpha(\cdot,\lambda,v_0,v_1,w)$ solves (\ref{eq:unpert:2}).
Moreover, if we take for given $(\gamma_0,\tilde{v}_0)$ 
the positive oriented orthonormal system $(v_0,v_1,w)$ defined by
\begin{align*}
v_0 &:= \lambda^{-1} \tilde{v}_0,\, v_1 := r\gamma_0+\sqrt{1-r^2}(v_0\cross \gamma_0),\,
w := (v_1\cross v_0)
\end{align*}
and $\lambda>0$ as above, then $\alpha(\cdot,\lambda,v_0,v_1,w)$ satisfies the initial conditions
\begin{align*}
\alpha(0,\lambda,v_0,v_1,w)=\gamma_0,\; \dot\alpha(0,\lambda,v_0,v_1,w)= \tilde{v}_0.  
\end{align*}
Since we are only interested in solutions in $H^{2,2}(S^1,S^2)$ we get an extra condition
on $\lambda$, i.e. the $1$-periodicity leads to
\begin{align*}
\lambda \in 2\pi \zz r.  
\end{align*}
Hence the simple solutions in $H^{2,2}(S^1,S^2)$ of (\ref{eq:unpert:1}) are
given by
\begin{align*}
\mathcal{Z} := \{\alpha(\cdot, 2\pi|r|,v_0,v_1,w ) \where &\{v_0,v_1,w\} \text{ is a}\\ 
&\text{positive orthonormal system in }\rz^3\}.  
\end{align*}
$S0(3)$ acts on solutions: 
if $\gamma$ solves (\ref{eq:unpert:1}) so does $A\circ\gamma$
for any $A\in SO(3)$. We have
\begin{align*}
A\circ \alpha(\cdot, 2\pi|r|,v_0,v_1,w ) = \alpha(\cdot, 2\pi|r|,A(v_0),A(v_1),A(w)),  
\end{align*}
and the set of solutions is parametrized by $SO(3)$.
It is easy to see, that
\begin{align*}
\alpha(\cdot, 2\pi|r|,v_0,v_1,w ) = \theta*\alpha(\cdot, 2\pi|r|,v_0',v_1',w')
\end{align*}
for some $\theta \in S^1$ if and only if $w=w'$. 
Consequently the set of critical orbits is parametrized by $w \in S^2$.
In the sequel we fix $k_0>0$ and $r>0$.\\
To compute the kernel of $D_{g_0}X_{k_0,g_0}\eval_{\alpha}$ at  
$\alpha=\alpha(\cdot,2\pi r,v_0,v_1,w)$ for some fixed system $(v_0,v_1,w)$
we note that for $V \in  T_{\alpha}H^{2,2}(S^1,S^2)$
\begin{align*}
R_{g_0}(V,\dot\alpha)\dot\alpha = V|\dot \alpha|^{2}-\langle V,\dot\alpha\rangle \dot\alpha  
\end{align*}
and hence by \eqref{eq:dg_x_g_formula}
\begin{align}
\label{eq:4}
D_{g_0}X_{k_0,g_0}\eval_{\alpha}(V) &= (-D_{g_0,t}^2+1)^{-1}
\big(-D_{t,g_0}^2V - V|\dot \alpha|^{2}+\langle V,\dot\alpha\rangle \dot\alpha \notag \\
&\quad + |\dot\alpha|^{-1}\langle D_{t,g_0}V,\dot\alpha\rangle k_0 (\alpha\cross \dot\alpha)
+|\dot\alpha|  k_0 (\alpha \cross D_{t,g_0}V)\big). 
\end{align}
Due to the geometric origin of equation (\ref{eq:unpert:1}) we deduce that
\begin{align*}
W_1(t,v_0,v_1,w) &:=  \dot\alpha= 2\pi r(-\sin(2\pi t)v_1 + \cos(2\pi t)v_0),\\
W_1(0,v_0,v_1,w) &=2\pi r v_0,\, D_{t,g_0} W_1(0,v_0,v_1,w)= 
-4\pi^2r^3k_0(k_0v_1-w),\\
W_0(t,v_0,v_1,w) &:= t\dot\alpha,\\
W_0(0,v_0,v_1,w)&=0 ,\, 
D_{t,g_0} W_0(0,v_0,v_1,w)= 2\pi r v_0,  
\end{align*}
solve the equation
\begin{align}
\label{eq:10}
0 =&-D_{t,g_0}^2W - W|\dot \alpha|^{2}+\langle W,\dot\alpha\rangle \dot\alpha \notag\\ 
&+ |\dot\alpha|^{-1}\langle D_{t,g_0}W,\dot\alpha\rangle k_0 (\alpha\cross \dot\alpha)
+|\dot\alpha|  k_0 (\alpha \cross D_{t,g_0}W). 
\end{align}
The vector-field
$W_1$ corresponds to invariance with respect to the $S^1$-action, 
$\theta \mapsto \alpha(\cdot+\theta)$,
and $W_0$ stems from the change of parameter $s \mapsto \alpha(\cdot s)$. The
$SO(3)$ invariance leads to two additional vector-fields in the kernel 
of $D_{g_0}X_{k_0,g_0}\eval_{\alpha}$, i.e.
we let
\begin{align*}
w_{1,s}:= \cos(s)w+sin(s)v_1,\, v_0=v_0,\\
v_{1,s}=v_0\cross w_{1,s}=\cos(s)v_1-\sin(s)w,\\
w_{2,s}:= \cos(s)w+sin(s)v_0,\, v_1=v_1,\\
v_{0,s}=w_{2,s}\cross v_1=\cos(s)v_0-\sin(s)w  
\end{align*}
and get
\begin{align}
\label{eq:formula_w2_w3}
W_2(t,v_0,v_1,w) &:= \frac{d}{ds}\big(\alpha(\cdot,2\pi r,v_0,v_{1,s},w_{1,s})\big\eval_{s=0}
= r k_0 v_1-r\cos(2\pi t)w, \notag \\
W_2(0,v_0,v_1,w) &= \sqrt{1-r^2} v_1-rw,\, D_{t,g_0} W_2(0,v_0,v_1,w)= 0,\notag\\
W_3(t,v_0,v_1,w) &:= \frac{d}{ds}\big(\alpha(\cdot,2\pi r,v_{0,s},v_{1},w_{2,s})\big\eval_{s=0}
= rk_0 v_0-r\sin(2\pi t)w,\notag\\
W_3(0,v_0,v_1,w) &= r k_0 v_0,\, D_{t,g_0} W_3(0,v_0,v_1,w)= 2\pi r^{3}(k_0v_0-w).  
\end{align}
We will omit the dependence of $W_i$ on $(v_0,v_1,w)$, if there is no possibility of confusion.
Since the initial values of $W_0,\dots, W_3$ are linearly independent in 
$\big(T_{\alpha(0)}S^2\big)^2$, any solution of \eqref{eq:10} is a linear combination
of $W_0,\dots, W_3$.
As only $W_1,\dots, W_3$ are periodic, we obtain
\begin{align}
\label{eq:kernel_DX_0}
\text{kernel}(D_{g_0}X_{k_0,g_0}\eval_{\alpha})=\langle W_1,\,W_2,\,W_3 \rangle
= T_{\alpha}\mathcal{Z}.  
\end{align}
To find the image of $D_{g_0}X_{k_0,g_0}\eval_{\alpha}$ we note 
that the moving frame
$\{\dot\alpha,\alpha\cross\dot\alpha\}$
is an orthogonal system in $T_{\alpha}S^{2}$ for any $t\in S^{1}$. 
Thus any $V \in T_\alpha H^{2,2}(S^1,S^2)$ may be written
as
\begin{align*}
V= \lambda_1 \dot\alpha + \lambda_2 (\alpha\cross \dot\alpha)   
\end{align*}
for some functions $\lambda_1,\,\lambda_2 \in H^{2,2}(S^{1},\rz)$.
Using the fact that
\begin{align*}
D_{t,g_0} \dot\alpha &= |\dot \alpha| k_0(\alpha \cross \dot\alpha) \text{ and }
D_{t,g_0} (\alpha \cross \dot\alpha) = -|\dot \alpha| k_0 \dot\alpha,    
\end{align*}
we may express $D_{t,g_0} V$ and $(D_{t,g_0})^{2} V$ in terms of $\lambda_1$ and $\lambda_2$. 
This leads to
\begin{align}
\label{eq:d_x_g0_alpha}
D_{g_0}X_{k_0,g_0}\eval_{\alpha}(V) &= (-D_{t,g_0}^2+1)^{-1}
\big((-\lambda_1''+2\pi\sqrt{1-r^{2}} \lambda_2') \dot\alpha \notag\\
&\qquad +(-\lambda_2''-(2\pi)^{2}\lambda_2)(\alpha\cross \dot\alpha)\big). 
\end{align}
Concerning $W_1,\dots, W_3$ and $W_g$ we find 
\begin{align}
W_1(t) &= \dot\alpha(t),\notag \\
W_2(t) &= -\frac{1}{2\pi r}\big(\sqrt{1-r^{2}} \sin(2\pi t)\dot\alpha(t) + \cos(2\pi t)(\alpha\cross
\dot\alpha)\big),\notag \\
W_3(t) &= -\frac{1}{2\pi r}\big(-\sqrt{1-r^{2}} \cos(2\pi t)\dot\alpha(t) + \sin(2\pi t)(\alpha\cross
\dot\alpha)\big)\notag \\
\label{eq:8}
W_{g_0}(\alpha) &= (1+|\dot\alpha|^{2}k_0^{2})^{-1}\dot\alpha =  (1+|\dot\alpha|^{2}k_0^{2})^{-1} W_1(\alpha). 
\end{align}

\begin{lemma}
\label{l:kernel_not_in_range}
For any solution $\alpha$ of the unperturbed problem there holds
\begin{align*}
\{0\}&=\langle W_1(\alpha),W_2(\alpha),W_3(\alpha)\rangle \cap R\big(D_{g_0}X_{k_0,g_0}\eval_{\alpha}\big) ,\\
\langle W_1(\alpha)\rangle^{\perp}&=  
\langle W_2(\alpha),W_3(\alpha)\rangle \oplus R\big(D_{g_0}X_{k_0,g_0}\eval_{\alpha}\big)
\end{align*}
\end{lemma}
\begin{proof}
We omit the dependence of $W_i$ on $\alpha$.
For $\lambda_1,\lambda_2 \in H^{2,2}(S^{1},\rz)$ we have
\begin{align*}
 (-D_{t,g_0}^2+1)&\big(\lambda_1 \dot\alpha + \lambda_2(\alpha \cross \dot\alpha)\big)\\
&= \big(-\lambda_1''+4\pi\sqrt{1-r^{2}} \lambda_2' +(4\pi^{2}(1-r^{2})+1) \lambda_1\big)\dot\alpha\\
&\quad +\big(-\lambda_2''-4\pi\sqrt{1-r^{2}} \lambda_1' +(4\pi^{2}(1-r^{2})+1) \lambda_2\big)\alpha \cross \dot\alpha
\end{align*}
Hence we get by direct calculations
\begin{align}
(-D_{t,g_0}^2+1)(W_1) &= (4\pi^{2}(1+r^{2}+1) \dot\alpha,\notag \\
(-D_{t,g_0}^2+1)(-2\pi r W_2) &= \sqrt{1-r^{2}}(-4\pi^{2}r^{2} +1) \sin(2\pi t) \dot\alpha \notag \\
\label{eq:D2_1_w2}
&\quad + (4\pi^{2} r^{2}+1)\cos(2\pi t) (\alpha \cross \dot\alpha),\\
(-D_{t,g_0}^2+1)(-2\pi r W_3) &= -\sqrt{1-r^{2}}(-4\pi^{2}r^{2} +1) \cos(2\pi t) \dot\alpha\notag\\
\label{eq:D2_1_w3}
&\quad + (4\pi^{2} r^{2}+1)\sin(2\pi t) (\alpha \cross \dot\alpha).
\end{align}
Consequently, by (\ref{eq:dx_to_wg_perp}) and \eqref{eq:8}
the vector $W_1$ is orthogonal  
to $\langle W_2,W_3\rangle$ and to $R\big(D_{g_0}X_{k_0,g_0}\eval_{\alpha}\big)$
in $T_\alpha H^{2,2}(S^1,S^2)$.
As in $L^{2}(S^{1},\rz)$
\begin{align*}
\lambda_2''+(2\pi)^{2}\lambda_2 \perp_{L^2} \langle \cos(2\pi t),\sin(2\pi t)\rangle,\;
\langle \lambda_1'',\lambda_2'\rangle  \perp_{L^2} {\rm const},
\end{align*}
we get
\begin{align*}
\{0\} &=  (-D_{t,g_0}^2+1)\big(\langle W_1,W_2,W_3\rangle\big)\\
&\qquad \cap 
(-D_{t,g_0}^2+1)D_{g_0}X_{k_0,g_0}\eval_{\alpha}(T_{\alpha}H^{2,2}(S^{1},S^{2}))   
\end{align*}
and the claim follows for $D_{g_0}X_{k_0,g_0}\eval_{\alpha}$ is a Fredholm operator of index $0$.
\end{proof}
To analyze the image of $D_{g_0}X_{k_0,g_0}$ we see for $\alpha \in \mathcal{Z}$
\begin{align}
\label{eq:decomp_range_dx}
R(D_{g_0}X_{k_0,g_0}\eval_{\alpha}) &= 
\big\{
(-D_{t,g_0}^2+1)^{-1}\big(
(-\lambda_1''+2\pi\sqrt{1-r^2}\lambda_2')\dot{\alpha} \notag\\
&\qquad -(\lambda_2''+(2\pi)^2\lambda_2)(\alpha \cross \dot{\alpha})\big)
\where \lambda_1,\,\lambda_2 \in H^{2,2}(S^1,\rz)
\big\}\notag \\
&=
\big\{
(-D_{t,g_0}^2+1)^{-1}\big(
\lambda_1\dot{\alpha}+\lambda_2(\alpha \cross \dot{\alpha})\big)
\where
\lambda_{1,2} \text{ in} \notag \\
&\qquad  L^2(S^1,\rz),\, \lambda_1 \perp_{L^2} 1,\,
\lambda_2 \perp_{L^2} \langle \cos(2\pi t),\sin(2\pi t)\rangle
\big\} \notag \\
&= \langle (\alpha \cross \dot{\alpha}) \rangle \oplus E_{+},
\end{align}
where $E_{+}$ is given by
\begin{align*}
E_{+} &=
\big\{
(-D_{t,g_0}^2+1)^{-1}\big(
\lambda_1\dot{\alpha}+\lambda_2(\alpha \cross \dot{\alpha})\big)
\where\\
&\qquad \lambda_1,\,\lambda_2 \in L^2(S^1,\rz),\, \lambda_1 \perp_{L^2}1,\,
\lambda_2 \perp_{L^2} \langle 1, \cos(2\pi t),\sin(2\pi t)\rangle
\big\}   
\end{align*}
We have for $V=\lambda_1\dot{\alpha}+\lambda_2 (\alpha \cross \dot{\alpha})$ in $T_\alpha H^{2,2}(S^{1},S^{2})$ 
\begin{align*}
D_{g_0}X_{k_0,g_0}\eval_{\alpha}(V) \in E_{+} \aequi \lambda_2 \perp_{L^{2}} 1
\aequi V \perp_{L^{2}} (\alpha \cross \dot{\alpha}).  
\end{align*}
We fix $V=(-D_{t,g_0}^{2}+1)^{-1}(\lambda_1\dot{\alpha}+\lambda_2 (\alpha \cross \dot{\alpha})) \in E_{+}$. Then
\begin{align*}
\int_{S^{1}}V (\alpha \cross \dot{\alpha}) &=
\int_{S^{1}} (-D_{t,g_0}^{2}+1)^{-1}(\lambda_1\dot{\alpha}+\lambda_2 (\alpha \cross \dot{\alpha})) (\alpha \cross
\dot{\alpha})\\
&= \int_{S^{1}} (\lambda_1\dot{\alpha}+\lambda_2 (\alpha \cross \dot{\alpha})) 
(-D_{t,g_0}^{2}+1)^{-1}(\alpha \cross \dot{\alpha})\\
&= (4\pi^{2}(1-r^{2})+1)^{-1}\int_{S^{1}} (\lambda_1\dot{\alpha}+\lambda_2 (\alpha \cross \dot{\alpha})) 
(\alpha \cross \dot{\alpha})=0.
\end{align*}
Consequently, $D_{g_0}X_{k_0,g_0}\eval_{\alpha}(E_{+})=E_{+}$.\\ 
Since $E_{+}$ is $L^2$-orthogonal to $\alpha \cross \dot{\alpha}$
and $\dot{\alpha}$, we may write
\begin{align*}
V = (\nu_1+f_1)\dot{\alpha}+(\nu_2+f_2)(\alpha \cross \dot{\alpha}),
\end{align*}
with
$\nu_{1,2} \perp_{L^2} \langle 1, \sin(2\pi \cdot),\cos(2\pi \cdot)\rangle$
and
$f_{1,2} \in \langle \sin(2\pi \cdot),\cos(2\pi \cdot)\rangle$.  
Then
\begin{align}
\label{eq:2}
\langle (-D_{t,g_0}^{2}&+1) D_{g_0}X_{k_0,g_0}\eval_{\alpha}(V),V\rangle_{L^{2}} \notag\\
 &=
\int_{S^{1}}
(\nu_1')^{2}-2\pi \sqrt{1-r^{2}} \nu_1'\nu_2
+(\nu_2')^{2}-4\pi^{2}(\nu_2)^{2} \notag\\
&\qquad +(f_1')^{2}-2\pi \sqrt{1-r^{2}} f_1'f_2.
\end{align}
For $\nu_2 \perp \langle 1,\, \cos(2\pi\cdot),\sin(2\pi\cdot)\rangle$ we have
\begin{align*}
\int_{S^{1}}(\nu_2')^{2}-4\pi^{2}(\nu_2)^{2}\ge \int_{S^{1}}16\pi^{2}(\nu_2)^{2},  
\end{align*}
hence
\begin{align*}
\int_{S^{1}}
(\nu_1')^{2}-2\pi \sqrt{1-r^{2}} \nu_1'\nu_2
+(\nu_2')^{2}-4\pi^{2}(\nu_2)^{2} 
&\ge \frac34 (\nu_1')^{2}+12 \pi^{2}(\nu_2)^{2}.
\end{align*}
Concerning the remaining term in \eqref{eq:2} we note that
as $(-D_{t,g_0}^{2}+1)$ maps 
\begin{align*}
\big\{\lambda_1 \dot{\alpha}+ \lambda_2 (\alpha \cross \dot{\alpha}) \where
\lambda_1,\lambda_2 \in \langle \sin(2\pi \cdot),\cos(2\pi \cdot)\rangle 
\big\}
\end{align*}
into itself and $V \in E_{+}$ there holds 
\begin{align*}
f_1 \dot{\alpha}+f_2(\alpha \cross \dot{\alpha})
\in
(-D_{t,g_0}^{2}+1)^{-1}\big\langle
\big(\cos(2\pi \cdot) \dot\alpha\big),
\big(\sin(2\pi \cdot) \dot\alpha\big)
\big\rangle. 
\end{align*}
Hence, by explicit computations there are $x,y \in \rz$ satisfying
\begin{align*}
f_1(t) &= x \cos(2\pi t)+ y \sin(2\pi t),\\
f_2(t) &= \frac{8\pi^2\sqrt{1-r^2}}{4\pi^2(2-r^2)+1}
\big(
y\cos(2\pi t)-x\sin(2\pi t)
\big).  
\end{align*}
This gives
\begin{align*}
\int_{S^{1}} (f_1')^2-2\pi \sqrt{1+r^2}f_1'f_2
 = \frac{2\pi^2(1+4\pi^2r^2)}{4\pi^2(2-r^2)+1}(x^2+y^2).
\end{align*}
This shows that
\begin{align*}
\langle (-D_{t,g_0}^{2}+1) D_{g_0}X_{k_0,g_0}\eval_{\alpha}(V),V\rangle_{L^{2}} >0
\text{ for all } V \in E_{+}\setminus \{0\},  
\end{align*}
and the homotopy
\begin{align*}
[0,1] \ni s \mapsto (1-s) \big(D_{g_0}X_{k_0,g_0}\eval_{\alpha}\big)\eval_{E_{+}} +s\; id\eval_{E_{+}}  
\end{align*}
is admissible. We use the decomposition in \eqref{eq:decomp_range_dx} and
\begin{align*}
D_{g_0}X_{k_0,g_0}\eval_{\alpha}(\alpha \cross \dot{\alpha}) = -\frac{4 \pi^2}{4\pi^2(1-r^2)+1} (\alpha \cross \dot{\alpha})   
\end{align*}
to see that   
\begin{align}
\label{eq:sgn_dx_g0_range}
\sgn \big(D_{g_0}X_{k_0,g_0}\eval_{\alpha}\big)\eval_{R(D_{g_0}X_{k_0,g_0}\eval_{\alpha})} =-1.
\end{align}

\subsection{The finite dimensional reduction}
\label{subsec:finite_reduction}
We fix $\alpha_0 \in \mathcal{Z}$ and a para\-metrization $\phi$ of $\mathcal{Z}$, 
which maps an open neighborhood of 
$0$ in $T_{\alpha_0}\mathcal{Z}$
into $\mathcal{Z}$, such that
\begin{align*}
\phi(0)=\alpha_0 \text{ and } D\phi\eval_0=id.  
\end{align*}
As $\mathcal{Z}$ consists of smooth functions, 
$\mathcal{Z}$ is a sub-manifold of $H^{m,2}(S^{1},S^{2})$ for $1\le m<\infty$.
We define $\Phi$ from an open neighborhood $\mathcal{U}$
of $0$ in
\begin{align*}
T_{\alpha_0}H^{2,2}(S^{1},S^{2})=\langle W_1(\alpha_0),W_2(\alpha_0),W_3(\alpha_0)\rangle \oplus 
R(DX_{g_0,0}\eval_{\alpha_0})  
\end{align*}
to $H^{2,2}(S^{1},S^{2})$ by 
\begin{align*}
\Phi(W,U):= Exp_{\alpha_0,g_0}\big(Exp_{\alpha_0,g_0}^{-1}(\phi(W))+U\big).  
\end{align*}
Then $(\Phi,\mathcal{U})$ is a chart
of $H^{2,2}(S^{1},S^{2})$ around $\alpha_0$ such that
$\mathcal{U}$ is an open neighborhood of $0$ in $T_{\alpha_0}H^{2,2}(S^{1},S^{2})$, and 
\begin{align*}
\Phi(0)=\alpha_0,\, D\Phi\eval_0=id,\\ 
\Phi^{-1}\big(\mathcal{Z}\cap \Phi(\mathcal{U})\big)= \mathcal{U}\cap \langle W_1(\alpha_0),W_2(\alpha_0),W_3(\alpha_0)\rangle.
\end{align*}
From the properties of $Exp_{\alpha_0,g_0}$ the map $\Phi$ is a chart of of $H^{k,2}(S^{1},S^{2})$ around
$\alpha_0$ for any $1\le k \le 4$ and shrinking $\mathcal{U}$ we may assume that
(\ref{eq:def_trans_1})-(\ref{eq:def_trans_3}) continue to hold with $Exp_{\gamma,g}$ replaced by $\Phi$, i.e.
\begin{align}
\label{eq:def_trans_1_phi}
T_{\Phi(V)}H^{1,2}(S^{1},S^{2}) 
&= \langle \frac{d}{dt}\Phi(V)\rangle \oplus D\Phi\eval_V(\langle \dot\alpha_0\rangle^{\perp,H^{1,2}}),\\
\label{eq:def_trans_2_phi}
T_{\Phi(V)}H^{2,2}(S^{1},S^{2}) 
&= \langle W_{g_0}(\Phi(V))\rangle \oplus D\Phi\eval_V(\langle W_{g_0}(\alpha_0)\rangle^{\perp}),\\
\label{eq:def_trans_3_phi}
\text{Proj}_{\langle W_{g_0}(\Phi(V)\rangle^\perp}\circ D \Phi\eval_{V}&:\,
\langle W_{g_0}(\alpha_0)\rangle^\perp \xpfeil{\iso}{} \langle W_{g_0}(\Phi(V)\rangle^\perp,
\end{align}
and the norm of the projections in (\ref{eq:def_trans_1_phi}) and (\ref{eq:def_trans_2_phi})
as well as the norm of the map in (\ref{eq:def_trans_3_phi}) and its inverse are uniformly bounded with respect to $V$.
For $\alpha_0 \in \mathcal{Z}$ the vectors $W_1(\alpha_0)$ and $W_{g_0}(\alpha_0)$ are collinear
and we use $\langle W_1(\alpha_0)\rangle$ instead of $\langle W_{g_0}(\alpha_0)\rangle$
in the analysis of the unperturbed problem below.\\
As in (\ref{eq:def_psi}) we get 
a chart $\Psi$ for the bundle $SH^{2,2}(S^{1},S^{2})$ around $(\alpha_0,0)$,
\begin{align*}
\Psi: \mathcal{U}\times \mathcal{U}\cap \langle W_1(\alpha_0)\rangle^{\perp}
\to SH^{2,2}(S^{1},S^{2}),\\
\Psi(V,U):= \big(\Phi(V),Proj_{\langle W_{g_0}(\Phi(V))\rangle^{\perp}}\circ D\Phi\eval_V(U)\big).
\end{align*}
Analogous to (\ref{eq:def_x_gamma}) we define
\begin{align*}
X_{g_0,\eps}^{\Phi}:\: \mathcal{U}\cap \langle W_1(\alpha_0)\rangle^{\perp}
\to \langle W_1(\alpha_0)\rangle^{\perp}   
\end{align*}
by
\begin{align*}
X_{g_0,\eps}^{\Phi}(V) := Proj_2 \circ \Psi^{-1}\big(\Phi(V),X_{g_0,\eps}(\Phi(V))\big).  
\end{align*}
Replacing $Exp_{\gamma,g}$ by $\Phi$ it is easy to see that
Lemma \ref{l:nondegenerate} carries over to $X_{g_0,\eps}^{\Phi}$, i.e.
\begin{align}
\label{eq:nondegenerate_phi}
V \in \mathcal{U}\cap \langle W_1(\alpha_0)\rangle^{\perp}
\text{ is a (nondegenerate) zero of } X_{g_0,\eps}^{\Phi} \text{ if and only if } \notag \\
S^{1}*\Phi(V) \text{ is a (nondegenerate) critical orbit of } X_{g_0,\eps},
\end{align}
and
if $X_{g_0,\eps}^{\Phi}(V)=0$, then after shrinking $\mathcal{U}$
\begin{align}
\label{eq:d_x_phi}
DX_{g_0,\eps}^{\Phi}\eval_V = A_V^{-1}\circ DX_{g_0,\eps}\eval_{\Phi(V)} \circ D\Phi\eval_V,  
\end{align}
where the isomorphism $A_V:\: \langle W_1(\alpha_0)\rangle^{\perp} \to \langle W_{g_0}(\Phi(V))\rangle^{\perp}$ 
is given by
\begin{align*}
A_V = Proj_{\langle W_{g_0}(\Phi(V))\rangle^{\perp}}\circ D\Phi\eval_V.  
\end{align*} 
From Lemma \ref{l:kernel_not_in_range} we may assume
\begin{align*}
\mathcal{U}\cap \langle W_1(\alpha_0)\rangle^{\perp} = \mathcal{U}_1\times \mathcal{U}_2,   
\end{align*}
where $\mathcal{U}_1$ and $\mathcal{U}_2$ are open neighborhoods of $0$ 
in $\langle W_2(\alpha_0),W_3(\alpha_0)\rangle$
and $R\big(D_{g_0}X_{k_0,g_0}\eval_{\alpha_0}\big)$.  
We denote for $\alpha \in \mathcal{Z}$ by $P_2(\alpha)$ the projection onto 
$R(DX_{g_0,0}\eval_{\alpha})$ with respect to the decomposition
\begin{align*}
\langle W_1(\alpha)\rangle^{\perp}&=  
\langle W_2(\alpha),W_3(\alpha)\rangle \oplus R\big(D_{g_0}X_{k_0,g_0}\eval_{\alpha}\big),
\end{align*}
and by $P_1(\alpha)$ the projection onto $\langle W_2(\alpha),W_3(\alpha)\rangle$. 
Moreover, for $W \in \mathcal{U}_1$ we define for $i=1,2$
\begin{align*}
P_i^{\Phi}(W):= (A_W)^{-1}\circ P_i(\Phi(W)) \circ A_W.
\end{align*}
The projections $P_1^{\Phi}(W)$ and $P_2^{\Phi}(W)$ correspond to the decomposition
\begin{align}
\label{eq:decomp_phi_w}
\langle W_1(\alpha_0)\rangle^{\perp}&=  
\langle W_2(\alpha_0),W_3(\alpha_0)\rangle \oplus R\big(D_{g_0}X_{g_0,0}^{\Phi}\eval_{W}\big),
\end{align}
as we have for $W \in \mathcal{U}_1$
\begin{align*}
DX_{g_0,0}^{\Phi}\eval_W = A_W^{-1}\circ DX_{g_0,0}\eval_{\Phi(W)} \circ A_W.  
\end{align*}

\begin{lemma}
\label{l:implicit_function}
For $\alpha_0\in \mathcal{Z}$ after possibly shrinking $\mathcal{U}$
there are $\eps_0>0$ and 
\begin{align*}
U &\in C^2([-\eps_0,\eps_0]\times \mathcal{U}_1,\langle W_1(\alpha_0)\rangle^\perp),\\
R &\in C^2([-\eps_0,\eps_0]\times \mathcal{U}_1,\langle W_2(\alpha_0),W_2(\alpha_0)\rangle),
\end{align*}
such that for all $(\eps,W) \in [-\eps_0,\eps_0]\times \mathcal{U}_1$ 
\begin{align*}
R(\eps,W) &=X_{g_0,\eps}^{\Phi}(W+U(\eps,W)),\\
0&= P_1^\Phi(W)\circ U(\eps,W),\\
O(\eps)_{\eps \to 0} &=\|U(\eps,W)\|+\|D_W U(\eps,W)\|+\|R(\eps,W)\|+\|D_W R(\eps,W)\|,\\
R(\eps,W)&= \eps P_1^\Phi(W)\circ K_1^{\Phi}(W)+o(\eps)_{\eps \to 0},\\
U(\eps,W) &= -\eps (DX_{g_0,0}^\Phi\eval_{W})^{-1}\circ  P_2^\Phi(W)\circ K_1^{\Phi}(W)+o(\eps)_{\eps \to 0}.  
\end{align*}
Moreover, $U(\eps, W)$ and $R(\eps,W)$ are unique in the following sense:\\ 
If $(\eps,W,U,R)$ in $[-\eps_0,\eps_0]\times \mathcal{U}_1\times \mathcal{U}\cap\langle W_1(\alpha_0)\rangle^\perp
\times \mathcal{U}_1$ 
satisfies
\begin{align*}
X_{g_0,\eps}^{\Phi}(W+U)=R \text{ and }
P_1^\Phi(W)\big(U\big) =0,  
\end{align*}
then $U=U(\eps,W)$ and $R=R(\eps,W)$.
\end{lemma}
\begin{proof}
We define a $C^2$-function $H$
\begin{align*}
H:\, \rz \times \mathcal{U}_1 &\times 
\mathcal{U}\cap\langle W_1(\alpha_0)\rangle^\perp
\times \langle W_2(\alpha_0),W_3(\alpha_0)\rangle\\
&\to
\langle W_1(\alpha_0)\rangle^\perp\times \langle W_2(\alpha_0),W_3(\alpha_0)\rangle,
\end{align*}
by
\begin{align*}
H(\eps,W,U,R) := \big(X_{g_0,\eps}^{\Phi}(W+U)-R,P_1^\Phi(W)(U)\big). 
\end{align*}
We have in $\mathcal{L}(\langle W_1(\alpha_0)\rangle^\perp\times \langle W_2(\alpha_0),W_3(\alpha_0)\rangle)$
\begin{align*}
D_{(U,R)} H\eval_{(0,W,0,0)}=
\begin{pmatrix}
DX_{g_0,0}^\Phi\eval_{W} &-id\\
P_1^\Phi(W) &0  
\end{pmatrix} 
,
\end{align*}
where we used the fact that $X_{g_0,0}^\Phi(W)=0$ and (\ref{eq:d_x_phi}).
From (\ref{eq:kernel_DX_0}) and Lemma \ref{l:kernel_not_in_range} 
we see that $D_{(U,R)} H\eval_{(0,0,0,0)}$
is an isomorphism. By the implicit function theorem, after possibly shrinking $\mathcal{U}$,
we get $\eps_0>0$ and unique functions $U=U(\eps,W)$ and $R=R(\eps,W)$
such that $H(\eps,W,U(\eps,W),R(\eps,W))=0$ for all
$(\eps,W) \in [-\eps_0,\eps_0]\times \mathcal{U}_1$, and $D_{(U,R)} H\eval_{\eps,W,U,R}$ is
uniformly invertible for $(\eps,W,U,R)\in [-\eps_0,\eps_0]\times \mathcal{U}_1 \times \mathcal{U}_2$.
This yields the existence and uniqueness part of the claim.\\
The uniqueness implies $U(0,W)=0$ 
and $R(0,W)=0$ for all $W\in \mathcal{U}_1$.
As $U$ and $R$ are differentiable we find $U(\eps,W)=O(\eps)$ and $R(\eps,W)=O(\eps)$ as $\eps \to 0$.
Moreover, taking the derivative with respect to
$W$ we see
\begin{align*}
0 &= D_W H\eval_{(0,W,0,0)} + D_{(U,R)} H\eval_{(0,W,0,0)} \big(D_W U(0,W),D_W R(0,W)\big)^T
\end{align*}
Since $H(0,W,0,0)\equiv 0$ we have $D_W H\eval_{(0,W,0,0)}=0$, which implies
\begin{align*}
(D_W U(0,W),D_W R(0,W))=(0,0),  
\end{align*}
because
$D_{(U,R)} H\eval_{(0,W,0,0)}$ is invertible. This gives the desired estimate for $D_WU$ and
$D_W R$.\\
Moreover, taking the derivative with respect to $\eps$ at $(0,W,0,0)$ we see as above
\begin{align*}
0 &= D_\eps H\eval_{(0,W,0,0)} + D_{(U,R)} H\eval_{(0,W,0,0)} \big(D_\eps U(0,W),D_\eps R(0,W)\big)^T\\
&= (K_1^{\Phi}(W),0)
+\begin{pmatrix}
DX_{g_0,0}^\Phi\eval_{W} &-id\\
P_1^\Phi(W) &0  
\end{pmatrix}
\begin{pmatrix}
D_\eps U(0,W)\\D_\eps R(0,W)  
\end{pmatrix}
\end{align*}
Consequently,
\begin{align*}
D_\eps R(0,W) &= P_1^\Phi(W)\circ K_1^{\Phi}(W),\\
D_\eps U(0,W) &= -(DX_{g_0,0}^\Phi\eval_{W})^{-1}\circ  P_2^\Phi(W)\circ K_1^{\Phi}(W) 
\end{align*}
This yields the claim.
\end{proof}

\begin{lemma}
\label{l:expansion_eps_0}
Under the assumptions of Lemma \ref{l:implicit_function} we have
as $\eps \to 0$
\begin{align*}
X_{g_0,\eps}^{\Phi}(W+U(\eps,W))
&= \eps P_1^\Phi(W)\circ K_1^{\Phi}(W) + O(\eps^{2})_{\eps \to 0},  
\end{align*}
where $K_1^{\Phi}$ is the vector-field $K_1$ in the coordinates $\Phi$, i.e.
\begin{align*}
K_1^{\Phi} = X_{g_0,1}^\Phi-X_{g_0,0}^\Phi.  
\end{align*}
\end{lemma}
\begin{proof}
Since $U(\eps,W)=O(\eps)$ we find
\begin{align*}
X_{g_0,\eps}^{\Phi}&(W+U(\eps,W)) \\
&= P_1^\Phi(W)\circ X_{g_0,\eps}^{\Phi}(W+U(\eps,W))  \\
&= P_1^\Phi(W)\circ X_{g_0,0}^{\Phi}(W+U(\eps,W)) + \eps P_1^\Phi(W)\circ K_1^{\Phi}(W+U(\eps,W))\\
&= P_1^\Phi(W)\circ DX_{g_0,0}^{\Phi}\eval_{W}U(\eps,W) +
\eps P_1^\Phi(W)\circ K_1^{\Phi}(W) + O(\eps^2)\\
&= \eps P_1^\Phi(W)\circ K_1^{\Phi}(W) 
+ O(\eps^{2})_{\eps \to 0}.  
\end{align*}
\end{proof}

\begin{lemma}
\label{l:nondeg_zero}
Under the assumptions of Lemma \ref{l:implicit_function}
suppose $0$ is a nondegenerate zero of the vector-field $P_1^\Phi(\cdot)\circ K_1^\Phi(\cdot)$ on
$\mathcal{U}_1$,
in the sense that $P_1^\Phi(0)\circ K_1^\Phi(0)=0$
and
\begin{align*}
D_W(P_1^\Phi(\cdot)\circ K_1^\Phi(\cdot))\eval_0 
\in \mathcal{L}(\langle {W}_2(\alpha_0),{W}_3(\alpha_0)\rangle)  
\end{align*}
is an isomorphism. Then, after possibly shrinking $\eps_0$ and $\mathcal{U}$,
for any $0<|\eps|\le \eps_0$ there is a unique $W(\eps)\in \mathcal{U}_1$
such that
\begin{align*}
X_{g_0,\eps}^{\Phi}(W(\eps)+U(\eps,W(\eps)))=0,\\
W(\eps) \to 0 \text{ as } \eps \to 0.  
\end{align*}
Moreover, $V(\eps):= W(\eps)+U(\eps,W(\eps))$ is the only zero of
$X_{g_0,\eps}^{\Phi}$ in $\mathcal{U}\cap \langle W_1(\alpha_0)\rangle^{\perp}$
and is nondegenerate with
\begin{align*}
\sgn(D X_{g_0,\eps}^{\Phi}\eval_{V(\eps)}) 
&= -det(D_W(P_1^\Phi(\cdot)\circ K_1^\Phi(\cdot))\eval_0).  
\end{align*}
\end{lemma}
\begin{proof}
Using Lemma \ref{l:implicit_function} and the estimates for $U$ and $D_W U$ we find
\begin{align}
\label{eq:3}
D_W&\Big(X_{g_0,\eps}^{\Phi}(\cdot+U(\eps,\cdot))\Big)\eval_W \notag\\
&= D_W\Big(P_1^\Phi(\cdot)\circ X_{g_0,\eps}^{\Phi}(\cdot+U(\eps,\cdot))\Big) \notag \\
&= (D_WP_1^\Phi\eval_W) \circ X_{g_0,\eps}^{\Phi}(W+U(\eps,W)) \notag\\
&\quad +
P_1^\Phi(W)
\circ D X_{g_0,\eps}^{\Phi}\eval_{W+U(\eps,W)}\circ (Id+D_W U\eval_{(\eps,W)})
\notag \\
&= (D_WP_1^\Phi\eval_W) \circ 
\Big(\eps K_1^\Phi(W)+DX_{g_0,\eps}^{\Phi}\eval_W U(\eps,W) + O(\eps^2)\Big)\notag \\
&\quad +
P_1^\Phi(W) \circ
\big(\eps D K_1^\Phi\eval_{W}+
D^2 X_{g_0,0}^{\Phi}\eval_{W}U(\eps,W)
+O(\eps^2)\big)
\end{align}
Differentiating the identity for fixed $\eps$
\begin{align*}
P_1^{\Phi}(W)\circ D X_{g_0,0}^{\Phi}\eval_W U(\eps,W) \equiv 0   
\end{align*}
with respect to $W$ we obtain
\begin{align}
\label{eq:7}
0&=(D_W P_1^{\Phi}\eval_W)\circ D X_{g_0,0}^{\Phi}\eval_W U(\eps,W)\notag \\
&\quad + P_1^{\Phi}(W)\circ \Big(D^2 X_{g_0,0}^{\Phi}\eval_W U(\eps,W)+
D X_{g_0,0}^{\Phi}\eval_W \circ  D_W U\eval_{(\eps,W)}\Big).
\end{align}
Since $P_1^{\Phi}(W)\circ D X_{g_0,0}^{\Phi}\eval_W \equiv 0$, combining (\ref{eq:3}) and (\ref{eq:7}) leads to
\begin{align}
\label{eq:D_W_P1_X}
D_W\Big(X_{g_0,\eps}^{\Phi}(\cdot+U(\eps,\cdot))\Big)\eval_W 
&= \eps D_W\big(P_1^{\Phi}(\cdot)\circ K_1^{\Phi}(\cdot)\big)\eval_W
+O(\eps^{2}).  
\end{align}
We define $F:[-\eps_0,\eps_0]\times \mathcal{U}_1 \to \langle {W}_2(\alpha_0),{W}_3(\alpha_0)\rangle$
by
\begin{align*}
F(\eps,W) := \eps^{-1} P_1^\Phi(W)\circ X_{g_0,\eps}^{\Phi}(W+U(\eps,W)).  
\end{align*}
Note that by Lemma \ref{l:expansion_eps_0} the function $F$ extends continuously to $\eps=0$.
By (\ref{eq:D_W_P1_X}) we have
\begin{align*}
D_W F\eval_{(\eps,W)} &= D_W\big(P_1^{\Phi}(\cdot)\circ K_1^{\Phi}(\cdot)\big)\eval_W
+O(\eps),  
\end{align*}
and $F$ is in $C^1$ with $D_W F\eval_{(0,0)}$ invertible. 
Consequently, by the implicit function theorem
after shrinking $\eps_0$ and $\mathcal{U}$
there is a unique $C^1$-function $W=W(\eps)$ such that $F(\eps,W(\eps))\equiv 0$ and
for $\eps\neq 0$
\begin{align*}
X_{g_0,\eps}^{\Phi}(W(\eps)+U(\eps,W(\eps))) \equiv 0.  
\end{align*}
Shrinking $\mathcal{U}$ we may assume that any $V \in \mathcal{U}\cap \langle W_1(\alpha_0)\rangle^{\perp}$
admits a unique decomposition $V=W_V+U_V$, where $U_V=P_2^{\Phi}(W_V)V$.   
From the construction in Lemma \ref{l:implicit_function} and the analysis above we see that
for $(\eps,V) \in [-\eps_0,\eps_0]\setminus\{0\}
\times \mathcal{U}\cap \langle W_1(\alpha_0)\rangle^{\perp}$ 
\begin{align*}
X_{g_0,\eps}^{\Phi}(V) = 0 &\aequi X_{g_0,\eps}^{\Phi}(W_V+U_V) = 0\\
&\aequi U_V=U(\eps,W_V) \text{ and }  X_{g_0,\eps}^{\Phi}(W_V+U(\eps,W_V)) = 0\\
&\aequi V=W(\eps)+U(\eps,W(\eps)).
\end{align*}
We use the decomposition in (\ref{eq:decomp_phi_w})
to compute the local degree of $X_{g_0,\eps}^{\Phi}$ in $V(\eps):=W(\eps)+U(\eps,W(\eps))$
as $\eps \to 0$. As $U(\eps,W)=O(\eps)$ we find
\begin{align}
\label{eq:expansion_d_x_phi}
D X_{g_0,\eps}^{\Phi}\eval_{V(\eps)} &= 
D X_{g_0,0}^{\Phi}\eval_{W(\eps)}+ D^{2} X_{g_0,0}^{\Phi}\eval_{W(\eps)}U(\eps,W(\eps))
\notag \\
&\qquad + \eps D K_1^{\Phi}\eval_{W(\eps)} + O(\eps^{2}) 
\end{align}
Differentiating for fixed $\tilde{W}\in \langle W_2(\alpha_0),W_3(\alpha_0)\rangle$
the identity
$$
D X_{g_0,0}^{\Phi}\eval_{W}\tilde{W} \equiv 0,  
$$
we obtain
$
D^{2} X_{g_0,0}^{\Phi}\eval_{W}\tilde{W} \equiv 0  
$
and thus by (\ref{eq:expansion_d_x_phi})
\begin{align*}
D X_{g_0,\eps}^{\Phi}\eval_{V(\eps)}\tilde{W} &= 
\big(\eps D K_1^{\Phi}\eval_{W(\eps)} + O(\eps^{2})\big)\tilde{W}. 
\end{align*}
For $\tilde{U}\in R\big(D_{g_0}X_{g_0,0}^{\Phi}\eval_{W}\big)$ we get from
(\ref{eq:expansion_d_x_phi})
\begin{align*}
D X_{g_0,\eps}^{\Phi}\eval_{V(\eps)}\tilde{U} &= 
\big(D X_{g_0,0}^{\Phi}\eval_{W(\eps)} + O(\eps)\big)\tilde{U}.  
\end{align*}
Consequently, with respect to the decomposition in (\ref{eq:decomp_phi_w})
\begin{align*}
D X_{g_0,\eps}^{\Phi}\eval_{V(\eps)} &=
\begin{pmatrix}
\eps P_1^{\Phi}(W(\eps)) \circ D K_1^{\Phi}\eval_{W(\eps)} & 0\\
0 &  D X_{g_0,0}^{\Phi}\eval_{W(\eps)}
\end{pmatrix}
\\
&\quad+
\begin{pmatrix}
O(\eps^{2}) & O(\eps)\\
O(\eps) & O(\eps)  
\end{pmatrix}.
\end{align*}
This shows that shrinking $\eps_0>0$ we may assume that $V(\eps)$ is a nondegenerate zero
of $X_{g_0,\eps}^{\Phi}$ for all $0<|\eps|\le \eps_0$ and by (\ref{eq:sgn_dx_g0_range})
\begin{align*}
\sgn(D X_{g_0,\eps}^{\Phi}\eval_{V(\eps)}) &= 
\det\big(D(P_1^{\Phi}(\cdot)\circ K_1^{\Phi}(\cdot))\eval_0\big)
\sgn(D X_{g_0,0}^{\Phi}\eval_{W(\eps)})\\
&= -\det\big(D(P_1^{\Phi}(\cdot)\circ K_1^{\Phi}(\cdot))\eval_0\big)  
\end{align*}
This finishes the proof.
\end{proof}
We consider $P_1(\cdot)\circ K_1(\cdot)$ as a vector field on $\mathcal{Z}$.
If $\alpha_0 \in \mathcal{Z}$ is a zero of $P_1(\cdot)\circ K_1(\cdot)$ then we obtain
due to the $S^1$ invariance that $S^1*\alpha_0 \subset \mathcal{Z}$ is
a zero orbit and
\begin{align*}
W_1(\alpha_0) \in \text{kernel}(D_{\mathcal{Z}}(P_1(\cdot)\circ K_1(\cdot))\eval_{\alpha_0}),\\
R(D_{\mathcal{Z}}(P_1(\cdot)\circ K_1(\cdot))\eval_{\alpha_0})\perp W_1(\alpha_0),
\end{align*}
where $D_{\mathcal{Z}}$ denotes the covariant
derivative on $\mathcal{Z}$. In the sequel we will therefore consider
$D_{\mathcal{Z}}(P_1(\cdot)\circ K_1(\cdot))\eval_{\alpha_0}$ as a map
\begin{align*}
D_{\mathcal{Z}}(P_1(\cdot)\circ K_1(\cdot))\eval_{\alpha_0}:\:
\langle {W}_2(\alpha_0),{W}_3(\alpha_0)\rangle \to 
\langle {W}_2(\alpha_0),{W}_3(\alpha_0)\rangle.  
\end{align*}
\begin{lemma}
\label{l:nondeg_zero_gamma}
Under the assumptions of Lemma \ref{l:implicit_function}
suppose $\alpha_0$ is a nondegenerate zero of the vector field 
$P_1(\cdot)\circ K_1(\cdot)$ on $\mathcal{Z}$,
in the sense that $P_1(\alpha_0)\circ K_1(\alpha_0)=0$
and
\begin{align*}
D_{\mathcal{Z}}(P_1(\cdot)\circ K_1(\cdot))\eval_{\alpha_0}
\in \mathcal{L}(\langle {W}_2(\alpha_0),{W}_3(\alpha_0)\rangle)  
\end{align*}
is an isomorphism. Then 
for any $0<\eps<\eps_0$ there is $\gamma(\eps) \in \Phi(\mathcal{U})$ satisfying
\begin{align*}
X_{g_0,\eps}(\gamma(\eps))=0 \text{ and }
\gamma(\eps) \to \alpha_0 \text{ as } \eps \to 0.  
\end{align*}
Moreover, $S^1*\gamma(\eps)$ is the unique critical orbit of $X_{g_0,\eps}$ 
in $\Phi(\mathcal{U})$
and is nondegenerate with
\begin{align*}
\deg_{loc, S^{1}}(X_{g_0,\eps},S^1*\gamma(\eps)) &= 
-\det(D_{\mathcal{Z}}(P_1(\cdot)\circ K_1(\cdot))\eval_{\alpha_0}).  
\end{align*}
\end{lemma}
\begin{proof}
We note that as $P_1(\alpha_0)\circ K_1(\alpha_0)=0$
\begin{align*}
D_{\mathcal{Z}}(P_1(\cdot)\circ K_1(\cdot))\eval_{\alpha_0} 
=
D(P_1^\Phi(\cdot)\circ K_1^\Phi(\cdot))\eval_{0}.  
\end{align*}
Consequently, the assumptions of Lemma \ref{l:nondeg_zero} are satisfied and
we may define for $0<\eps<\eps_0$ the curve $\gamma(\eps)$ by
\begin{align*}
\gamma(\eps)&:= \Phi(V(\eps)) \in H^{2,2}(S^1,S^2). 
\end{align*}
From \eqref{eq:nondegenerate_phi} we infer that $\gamma(\eps)$ is the unique zero
of $X_{g_0,\eps}$ in $\Phi(\langle W_1\rangle^\perp \cap \mathcal{U})$ and $S^1*\gamma(\eps)$
is a nondegenerate critical orbit.
It is easy to see that the existence of a slice in Lemma \ref{l:slice_lemma} remains valid
if we replace $Exp_{\alpha_0,g_0}$ by $\Phi$. Consequently, $S^1*\gamma(\eps)$ is the
unique critical orbit of $X_{g_0,\eps}$ in $S^1*\Phi(\langle W_1\rangle^\perp \cap \mathcal{U})$,
which is an open neighborhood of $S^1*\alpha_0$ in $H^{2,2}(S^1,S^2)$.\\
We fix $0<\eps<\eps_0$ and consider for $s\in[0,1]$ the family of maps
\begin{align*}
Y_s := A_{V(\eps)}^{-1}\circ DX_{g_0,\eps}\eval_{\gamma(\eps)}
\circ \big((1-s)+s\text{Proj}_{\langle W_1(\gamma(\eps))\rangle^\perp}\big)D\Phi\eval_{V(\eps)}.  
\end{align*}
Since $DX_{g_0,\eps}\eval_{\gamma(\eps)}$ restricted to $\langle W_1(\gamma(\eps))\rangle^\perp$
is of the form $id-compact$, writing 
\begin{align*}
DX_{g_0,\eps}\eval_{\gamma(\eps)} = DX_{g_0,\eps}\eval_{\gamma(\eps)}
\circ \text{Proj}_{\langle W_1(\gamma(\eps))\rangle^\perp}+
DX_{g_0,\eps}\eval_{\gamma(\eps)}
\circ \text{Proj}_{\langle W_1(\gamma(\eps))\rangle},  
\end{align*}
we deduce that $Y_s=id-compact$ for all $s\in [0,1]$.
From Lemma \ref{l:nondeg_zero} we have that $Y_0$ is invertible and satisfies
\begin{align*}
Y_0= DX_{g_0,\eps}^\Phi\eval_{V(\eps)} \text{ and } 
\sgn(Y_0)=-\det(D_{\mathcal{Z}}(P_1(\cdot)\circ K_1(\cdot))\eval_{\alpha_0}).  
\end{align*}
As $DX_{g_0,\eps}^\Phi\eval_{V(\eps)}$ is invertible, the kernel of $DX_{g_0,\eps}\eval_{\gamma(\eps)}$
is given by $\langle \dot{\gamma}(\eps)\rangle$.
Since $\gamma(\eps)$ converges to $\alpha_0$ as $\eps \to 0$ and $\dot{\alpha}_0=W_1(\alpha_0)$ 
we get
\begin{align*}
\dot{\gamma}(\eps)= W_1(\gamma(\eps))+o(1)_{\eps \to 0}, 
\end{align*}
which implies together with (\ref{eq:def_trans_1_phi}) 
that $\langle \dot{\gamma}(\eps)\rangle$ is transversal to the image of
\begin{align*}
\big((1-s)+s\text{Proj}_{\langle W_1(\gamma(\eps))\rangle^\perp}\big)\circ D\Phi\eval_{V(\eps)}  
\end{align*}
for all $s\in[0,1]$.
Consequently, $Y_s$ remains invertible when $s$ moves from $0$ to $1$. Due to the homotopy invariance
we finally obtain
\begin{align*}
\sgn(Y_0) &=
-\sgn\big(\det(D_{\mathcal{Z}}(P_1(\cdot)\circ K_1(\cdot))\eval_{\alpha_0})\big)\\
&=  \sgn(Y_1)
= \sgn\Big(A_{V(\eps)}^{-1}\circ DX_{g_0,\eps}\eval_{\gamma(\eps)} \circ A_{V(\eps)}
\Big)\\
&=\sgn(DX_{g_0,\eps}\eval_{\gamma(\eps)}) = \deg_{loc, S^{1}}(X_{g_0,\eps},S^{1}*\gamma(\eps)).  
\end{align*}
This finishes the proof.
\end{proof}

\subsection{The $S^1$-degree of \eqref{eq:unpert:1}}
\label{subsec:degree_unperturbed}
We define the function $k_1$ by
\begin{align}
\label{eq:def_k1}
k_1(x):= \langle x,e_3\rangle \text{ for } x \in S^{2}=\rand B_1(0) \subset \rz^{3}.  
\end{align}
The corresponding vector-field $K_1$ on $H^{2,2}(S^{1},S^{2})$ is given by
\begin{align*}
K_1(\alpha) = (-D_{t,g_0}^{2}+1)^{-1}(|\dot \alpha|\langle \alpha,e_3\rangle (\alpha \cross \dot{\alpha})).  
\end{align*}
We note that for $\alpha = \alpha(\cdot,2\pi r,v_0,v_1,w) \in \mathcal{Z}$ we have
\begin{align*}
(-&D_{t,g_0}^{2}+1)(2\pi r)^{-1} K_1(\alpha)\\
 &= 
\big(
\sqrt{1-r^{2}} \langle w,e_3\rangle + r \cos(2\pi \cdot)\langle v_1,e_3\rangle+
r \sin(2\pi \cdot)\langle v_0,e_3\rangle
\big)(\alpha \cross \dot{\alpha})\\
&= \frac{-2\pi r^{2}}{4\pi^{2}r^{2}+1}(-D_{t,g_0}^{2}+1)
\big(\langle v_1,e_3\rangle  W_2(\alpha)
+ \langle v_0,e_3\rangle W_3(\alpha)\big)
+ h(\alpha),  
\end{align*}
where $(-D_{t,g_0}^{2}+1)^{-1}h(\alpha)$ is in the image of $D_{g_0}X_{k_0,g_0}\eval{\alpha}$ by
(\ref{eq:d_x_g0_alpha})-(\ref{eq:D2_1_w3}).
Hence, 
\begin{align*}
P_1(\alpha)\circ K_1(\alpha) =   
\frac{-4\pi^2 r^{3}}{4\pi^{2}r^{2}+1}\langle v_1,e_3\rangle W_2(\alpha)
+ \frac{-4\pi^2 r^{3}}{4\pi^{2}r^{2}+1}\langle v_0,e_3\rangle W_3(\alpha),
\end{align*}
and there are exactly two critical orbits of 
$P_1(\alpha)\circ K_1(\alpha)$ on $\mathcal{Z}$ given by
\begin{align*}
\{\alpha= \alpha(\cdot,2\pi r,v_0,v_1,w) \in \mathcal{Z}\where w=\pm e_3\}
= S^1*\alpha_+\cup S^1*\alpha_-,
\end{align*}
where
\begin{align*}
\alpha_+ = \alpha(\cdot,2\pi r,e_1,e_2,e_3) \text { and }
\alpha_- = \alpha(\cdot,2\pi r,-e_1,e_2,-e_3).  
\end{align*}
The curves $\alpha_\pm$ correspond to two parallels 
with respect to the north pole $e_3$ and curvature $k_0$.
Using the formulas for $W_2$ and $W_3$ in \eqref{eq:formula_w2_w3} we find with respect to
the basis $\{W_2(\alpha_\pm),W_3(\alpha_\pm)\}$
\begin{align*}
D(P_1(\cdot)\circ K_1(\cdot))\eval_{\alpha_\pm} =\frac{4\pi^2 r^{3}}{4\pi^{2}r^{2}+1}
\begin{pmatrix}
\pm 1 & 0\\
0& \pm 1  
\end{pmatrix}.
\end{align*}
Thus, we may apply Lemma \ref{l:nondeg_zero_gamma} and get two critical orbits
$\alpha_\pm(\eps)$ for $X_{g_0,\eps}$ converging to $\alpha_\pm$ as $\eps \to 0$.

\begin{lemma}
\label{l:euler_charak_unperturbed}
Let $M$ be the subset of $H^{2,2}(S^{1},S^{2})$ consisting of simple and regular curves. Then
$\chi_{S^1}(X_{k_0,g_0},M)=-2$.  
\end{lemma}
\begin{proof}
We choose $k_1=\langle \cdot,e_3\rangle$ as above.
From Lemmas \ref{l:implicit_function}-\ref{l:nondeg_zero} there are $\eps_0>0$
and an open neighborhood
$\mathcal{U}$ of $\mathcal{Z}$ such that for all $0<\eps<\eps_0$ the critical orbits
of $X_{g_0,\eps}$ in $\mathcal{U}$ are given exactly by $S^1*\alpha_\pm(\eps)$.
Indeed, suppose there are $\eps_n \to 0^+$ and
a sequence $(\alpha_n)$ of zeros of $X_{g_0,\eps_n}$ converging to
$\mathcal{Z}$ different from $S^1*\alpha_\pm(\eps_n)$. Up to a subsequence
\begin{align*}
\alpha_n \to \alpha_0 \in \mathcal{Z}
\end{align*}
as $n \to \infty$. For large $n$ we use the chart $\Phi$ around $\alpha_0$ 
as in Lemma \ref{l:implicit_function}. 
From the existence of a slice in Lemma \ref{l:slice_lemma}
we get a sequence $\theta_n \in \rz/\zz$ converging to $0$ such that
\begin{align*}
\theta_n*\alpha_n = \Phi(V_n) \text{ for some }V_n \in \langle W_1(\alpha_0)\rangle^\perp.   
\end{align*}
As in the proof of Lemma \ref{l:nondeg_zero} we may decompose
\begin{align*}
V_n= \Phi^{-1}(\theta_n*\alpha_n)= W_n+U_n,  
\end{align*}
where $W_n \in \langle W_2(\alpha_0), W_3(\alpha_0)\rangle$ 
and $U_n \in R(DX_{k_0,g_0}^\Phi\eval_{W_n})$.
From the uniqueness part of Lemma \ref{l:implicit_function}, as $X_{g_0,\eps_n}(W_n+U_n)=0$,
we get $U_n= U(\eps_n,W_n)$. By Lemma \ref{l:expansion_eps_0} we see that necessarily
$P_1(\alpha_0)\circ K_1(\alpha_0)=0$, such that $S^1*\alpha_0 \in \{S^1*\alpha_{\pm}\}$.
From Lemma \ref{l:nondeg_zero_gamma} we finally
deduce that $S^1*\alpha_n \in \{ S^1*\alpha_{\pm}(\eps_n)\}$, a contradiction.\\
From the definition of the $S^1$-equivariant Poincar\'{e}-Hopf index and 
the classification of the simple zeros of $X_{k_0,g_0}$
there holds for small $\eps>0$ 
\begin{align*}
\chi_{S^1}(X_{k_0,g_0},M) &= \chi_{S^1}(X_{k_0,g_0},\mathcal{U})
= \chi_{S^1}(X_{g_0,\eps},\mathcal{U}) =-2. 
\end{align*}
\end{proof}

\section{Apriori estimates}
\label{sec:apriori}
We fix a continuous family of metrics $\{g_t\where t\in [0,1]\}$ on $S^2$
and a continuous family of positive continuous function $\{k_t \where t\in [0,1]\}$ on $S^2$.
We let $X_{t}$ be the vector field on $H^{2,2}(S^1,S^2)$ defined by
\begin{align*}
X_t:= X_{k_t,g_t}.  
\end{align*}
We denote by $M\subset H^{2,2}(S^1,S^2)$ the set
\begin{align*}
M:= \{\gamma \in H^{2,2}(S^1,S^2) \where \gamma \text{ is simple and regular.}\}.
\end{align*}
We shall give sufficient conditions assuring that the set 
\begin{align*}
X^{-1}(0) := \{(\gamma,t) \in M\times [0,1] \where X_t(\gamma)=0\}   
\end{align*}
is compact in $M\times [0,1]$.
Fix $(\gamma,t)\in X^{-1}(0)$.
The Gauss-Bonnet formula yields
\begin{align*}
\int_{\gamma}k_t \, ds +\int_{\Omega_\gamma}K_{g_t} \,dg_t = 2\pi,
\end{align*}
where $\Omega_\gamma$ denotes the interior of $\gamma$ with respect to 
the normal $N_{g_t}$ and $K_{g_t}$ is the Gauss curvature of $(S^2,g_t)$.
To obtain a contradiction assume that there is $(\gamma_n,t_n)$ in $X^{-1}(0)$
such that $L(\gamma_n)\to 0$ as $n \to \infty$. Then the left hand side in the Gauss-Bonnet
formula, as $k_t$ and $K_{g_t}$ are uniformly bounded, tends to $0$, which is impossible.  
Consequently, the length $L(\gamma)$ of $\gamma$
satisfies
\begin{align}
\label{eq:length_bound}
c\le L(\gamma) \le \big(\inf\{k_t(x)\}\big)^{-1} 
\big(2\pi + \sup_{t\in [0,1]}\{(\sup K_{g_t}^-) vol(S^2,g_t)\}\big),  
\end{align}
for some positive constant $c=c(\{k_t\},\{g_t\})$ and $K_{g_t}^-:=-\min(K_{g_t},0)$.\\
Suppose $(\gamma_n,t_n)$ in $X^{-1}(0)$ converges to $(\gamma_0,t_0)$ in $H^{2,2}(S^1,S^2)$, such that
\begin{align*}
\gamma_0 \not\in M.  
\end{align*}
Then by \eqref{eq:length_bound} the curve $\gamma_0$ is non-constant and regular, 
hence there is $s_1\neq s_2$ in $\rz/\zz$ such that
$\gamma_0(s_1)=\gamma_0(s_2)$. As $\gamma_n$ are simple curves, parametrized proportional to arc-length
we see that $\dot\gamma_0(s_1)=\pm\dot\gamma_0(s_2)$.
If $\dot\gamma_0(s_1)=\dot\gamma_0(s_2)$ then by the unique solvability
of the initial value problem
\begin{align*}
\gamma_0(\cdot+(s_1-s_2))=\gamma_0(\cdot).  
\end{align*}
If $\dot\gamma_0(s_1)=-\dot\gamma_0(s_2)$ then
we write $\gamma$ close to $s_1$ and $s_2$
as a graph over the tangent direction $\dot\gamma_0(s_1)$ in normal coordinates
$Exp_{\gamma_0(s_1)}$.
By the maximum principle we find
\begin{align*}
\gamma_0(s_1+t)&= Exp_{\gamma_0(s_1),g}\big(t\dot\gamma_0(s_1)+a(t)N_g(\gamma_0(s_1))\big),\\
\gamma_0(s_2+t)&= Exp_{\gamma_0(s_1),g}\big(-t\dot\gamma_0(s_1)-b(t)N_g(\gamma_0(s_1))\big),
\end{align*}
where $a(t)$ and $b(t)$ are positive for $t\neq 0$.
Consequently, if $\dot\gamma_0(s_1)=-\dot\gamma_0(s_2)$ then
$\gamma_0$ touches itself at $\gamma_0(s_1)$, locally 
separated by the geodesic through $\gamma_0(s_1)$ with velocity $\dot\gamma_0(s_1)$. 
Thus, $\gamma_0$ is a $m$-fold covering for some $m \in \nz$ of a curve $\alpha$, which is almost simple
in the sense that $\alpha$ can only touch itself as described above.
Using stereographic coordinates $\mathcal{S}$ 
there is a point $p_0$ close to the curve $\gamma_0$, such that
the winding number of $\mathcal{S}(\gamma_0)$ around $\mathcal{S}(p_0)$ 
is $\pm m$. Since $\gamma_0$ is a limit
of simple curves, by the stability of the winding number, we deduce $m=1$.\\
We denote by $(\Omega_0,g)$ the interior of $\gamma_0$ considered as a Riemannian surface with
boundary of positive geodesic curvature. Fix a touching point $\gamma_0(s_1)=\gamma_0(s_2)$.
The point $\gamma_0(s_1)=\gamma_0(s_2)$ corresponds to two different boundary points
of $\Omega_0$. Denote by $\beta$ the curve of minimal length in $\Omega_0$ 
connecting the two boundary points.
From a regularity result for variational problems with constraints (see \cite{MR620261,MR0420406})
the minimizer $\beta$ is a $C^{1}$-curve. By the maximum principle $\beta$ cannot touch
the boundary of $\Omega_0$ and is therefore a $C^{2}$ geodesic in the interior of $\Omega_0$.
Moreover, as a minimizer, $\beta$ is stable and going back to $S^2$
the curve $\beta$ is a geodesic loop which is stable with respect to variations with fixed
end-points. Thus
\begin{align}
\label{eq:bound_inj}
inj(g_{t_0}) \le \frac12 L(\beta) < \frac14 L(\gamma_0).
\end{align}
This leads to
\begin{lemma}
\label{l:x_0_closed}
$X^{-1}(0)$ is compact in $M\times [0,1]$ under each of the following assumptions
\begin{align}
\label{eq:cond_inj_unif}
\inf_{(t,x)\in [0,1]\times S^2}\{k_t\} 
\ge \frac14 \sup_{t\in [0,1]}\Big(\big(inj(g_t)\big)^{-1}\big(2\pi+(\sup K_{g_t}^-)vol(S^2,g_t)\big)\big),\\
\label{eq:cond_K_pos_unif}
K_{g_t}>0 \; \forall t\in [0,1] \text{ and } \inf_{(t,x)\in [0,1]\times S^2}\{k_t\}
 \ge \frac12 \sup_{t\in [0,1]}\Big(\big(\sup K_{g_t}\big)^{\frac12}\big),\\ 
\label{eq:cond_K_pinch_unif}
K_{g_t}>0 \; \forall t\in [0,1] \text{ and }
\big(\sup K_{g_t}\big) <4 \big(\inf K_{g_t}\big) \text{ for all }t\in [0,1],
\end{align}
where $inj(g_t)$ denotes the injectivity radius of $(S^2,g_t)$.  
\end{lemma}
\begin{proof}
We first show that $X^{-1}(0)$ is closed under each of the above assumptions.
Suppose $(\gamma_n,t_n)\in X^{-1}(0)$ converges to some $(\gamma_0,t_0)$ in $H^{2,2}(S^1,S^2)$.
To obtain a contradiction assume $(\gamma_0,t_0) \notin X^{-1}(0)$, i.e.
$\gamma_0$ is not simple. Then by the above analysis $\gamma_0$ touches itself at some
point $\gamma_0(s_1)=\gamma_0(s_2)$ and there is a stable, nontrivial geodesic loop $\beta$,
which yields a bound from above on the injectivity radius in \eqref{eq:bound_inj}
by the length of $\gamma_0$. If $\gamma_0$ is too short this is impossible.
The estimate on the length of $\gamma_0$ in \eqref{eq:length_bound} leads to the contradiction
under the assumption \eqref{eq:cond_inj_unif}.
If $K_{t_0}>0$ then by \cite[Thm 2.6.9]{MR1330918}
\begin{align}
\label{eq:bound_of_inj}
inj(g_{t_0}) \ge \pi \big(\sup K_{t_0}\big)^{-\frac12},
\end{align}
and (\ref{eq:cond_K_pos_unif}) is a special case of (\ref{eq:cond_inj_unif}).\\
Moreover, by Bonnet-Meyer's theorem, as $\beta$ is a stable geodesic loop, its length
is bounded by 
\begin{align*}
L(\beta) \le \frac{\pi}{\sqrt{\inf{K_{t_0}}}},  
\end{align*}
which yields together with (\ref{eq:bound_of_inj}) the contradiction
assuming (\ref{eq:cond_K_pinch_unif}).\\
To deduce the compactness of $X^{-1}(0)$ we fix a sequence $(\gamma_n,t_n)$ in $X^{-1}(0)$.
By (\ref{eq:length_bound}) the length $L_{g_{t_n}}(\gamma_n)$ is uniformly bounded.
Since each $\gamma_n$ is para\-metrized proportional to arc-length, $(|\dot \gamma_n|_{g_{t_n}})$
is uniformly bounded. Using the equation (\ref{eq:1}) and standard elliptic regularity
$(\gamma_n)$ is bounded in $H^{4,2}(S^{1},S^{2})$. Hence we may choose a subsequence, which converges
in $H^{2,2}(S^{1},S^{2})$ and by the first part of the proof in $X^{-1}(0)$ under 
each of the above assumptions.
This yields the claim. 
\end{proof}
\begin{proof}[{Proof of Theorem \ref{thm_example}}]
We fix $k_0>0$, let $k_1 \in C^\infty(S^2,\rz)$ be given by \eqref{eq:def_k1},
and consider the metrics $g_t\equiv g_0$, the functions $k_t := k_0+tk_1$, and
the corresponding vector fields $X_t := X_{k_t,g_0}$.
The zeros of $X_0$ in $M$ are given by $\mathcal{Z}$, the manifold of solutions to
the unperturbed problem. The compactness of $X^{-1}(0)$ implies that the zeros
of $X_t$ in $M$ converge to $\mathcal{Z}$ as $t \to 0$.
From the proof of Lemma \ref{l:euler_charak_unperturbed} there are
exactly two critical orbits $S^1*\alpha_\pm(t)$ for $|t|>0$ 
small enough close to $\mathcal{Z}$ which are nondegenerate and
converge to the orbits of the parallels $\alpha(\cdot, 2\pi|r|,\pm e_1,e_2,\pm e_3)$
as $t \to 0$. Consequently, there are exactly two simple solutions
of \eqref{eq:1} with $g=g_0$ and $k=k_0+tk_1$ if $|t|>0$ is small enough. 
\end{proof}

\section{Existence results}
\label{sec:existence}
We give the proof of our main existence result.
\begin{proof}[{Proof of Theorem \ref{thm_existence}}]
We consider the family of metrics $\{g_t\where t\in [0,1]\}$ defined by
\begin{align*}
g_t := (1-t)g_0+tg.  
\end{align*}
Since $\{g_t\}$ is a compact family of metrics, there is a constant $k_0>0$ such that
\begin{align*}
k_0 >
\frac14 \sup_{t\in [0,1]}\Big(\big(inj(g_t)\big)^{-1}\big(2\pi+(\sup K^-_{g_t})vol(S^2,g_t)\big)\big).  
\end{align*}
We denote by $M$ the set of simple regular curves in $H^{2,2}(S^1,S^2)$.
From condition \eqref{eq:cond_inj_unif} in Lemma \ref{l:x_0_closed} the homotopy 
\begin{align*}
[0,1]\ni t\mapsto X_{k_0,g_t}  
\end{align*}
is $(M,g_t,S^1)$-admissible and hence from Lemma \ref{l:degree_homotopy} and 
Lemma \ref{l:euler_charak_unperturbed}
\begin{align*}
-2 = \chi_{S^1}(X_{k_0,g_0},M)=\chi_{S^1}(X_{k_0,g},M).   
\end{align*}
We define the family of functions $\{k_t\where t\in [0,1]\}$ by
\begin{align*}
k_t := (1-t)k_0+tk
\end{align*}
and consider the homotopy
\begin{align*}
[0,1]\ni t\mapsto X_{k_t,g}.  
\end{align*}
Under each of the above assumptions we may apply Lemma \ref{l:x_0_closed}
to deduce that the homotopy is $(M,g,S^1)$-admissible, and thus
\begin{align*}
-2 = \chi_{S^1}(X_{k_0,g},M) = \chi_{S^1}(X_{k,g},M).  
\end{align*}
Since the local degree of an isolated critical orbit
is larger than $-1$ by Lemma \ref{l:degree_poincare}, there
are at least two simple solutions to \eqref{eq:1}. 
\end{proof}

\bibliographystyle{plain}
\bibliography{geodesic_curves}

\def\cprime{$'$} \def\cprime{$'$} \def\cprime{$'$}
\begin{thebibliography}{10}

\bibitem{MR515141}
Ralph Abraham and Jerrold~E. Marsden.
\newblock {\em Foundations of mechanics}.
\newblock Benjamin/Cummings Publishing Co. Inc. Advanced Book Program, Reading,
  Mass., 1978.
\newblock Second edition, revised and enlarged, With the assistance of Tudor
  Ra{\c{t}}iu and Richard Cushman.

\bibitem{MR620261}
Ralph Alexander and S.~Alexander.
\newblock Geodesics in {R}iemannian manifolds-with-boundary.
\newblock {\em Indiana Univ. Math. J.}, 30(4):481--488, 1981.

\bibitem{MR0420406}
F.~J. Almgren, Jr.
\newblock Existence and regularity almost everywhere of solutions to elliptic
  variational problems with constraints.
\newblock {\em Mem. Amer. Math. Soc.}, 4(165):viii+199, 1976.

\bibitem{MR890489}
V.~I. Arnold.
\newblock The first steps of symplectic topology.
\newblock {\em Uspekhi Mat. Nauk}, 41(6(252)):3--18, 229, 1986.

\bibitem{MR2078115}
Vladimir~I. Arnold.
\newblock {\em Arnold's problems}.
\newblock Springer-Verlag, Berlin, 2004.
\newblock Translated and revised edition of the 2000 Russian original, With a
  preface by V. Philippov, A. Yakivchik and M. Peters.

\bibitem{MR2036336}
Gonzalo Contreras, Leonardo Macarini, and Gabriel~P. Paternain.
\newblock Periodic orbits for exact magnetic flows on surfaces.
\newblock {\em Int. Math. Res. Not.}, (8):361--387, 2004.

\bibitem{MR2282341}
M.~C. Crabb.
\newblock The fibrewise {L}eray-{S}chauder index.
\newblock {\em J. Fixed Point Theory Appl.}, 1(1):3--30, 2007.

\bibitem{MR817033}
E.~N. Dancer.
\newblock A new degree for {$S\sp 1$}-invariant gradient mappings and
  applications.
\newblock {\em Ann. Inst. H. Poincar\'e Anal. Non Lin\'eaire}, 2(5):329--370,
  1985.

\bibitem{MR0277009}
K.~D. Elworthy and A.~J. Tromba.
\newblock Degree theory on {B}anach manifolds.
\newblock In {\em Nonlinear Functional Analysis (Proc. Sympos. Pure Math., Vol.
  XVIII, Part 1, Chicago, Ill., 1968)}, pages 86--94. Amer. Math. Soc.,
  Providence, R.I., 1970.

\bibitem{MR0209600}
F.~Brock Fuller.
\newblock An index of fixed point type for periodic orbits.
\newblock {\em Amer. J. Math.}, 89:133--148, 1967.

\bibitem{MR902290}
V.~L. Ginzburg.
\newblock New generalizations of {P}oincar\'e's geometric theorem.
\newblock {\em Funktsional. Anal. i Prilozhen.}, 21(2):16--22, 96, 1987.

\bibitem{MR1432462}
Viktor~L. Ginzburg.
\newblock On closed trajectories of a charge in a magnetic field. {A}n
  application of symplectic geometry.
\newblock In {\em Contact and symplectic geometry (Cambridge, 1994)}, volume~8
  of {\em Publ. Newton Inst.}, pages 131--148. Cambridge Univ. Press,
  Cambridge, 1996.

\bibitem{MR1417851}
Viktor~L. Ginzburg.
\newblock On the existence and non-existence of closed trajectories for some
  {H}amiltonian flows.
\newblock {\em Math. Z.}, 223(3):397--409, 1996.

\bibitem{MR2060021}
Viktor~L. Ginzburg and Ba{\c{s}}ak~Z. G{\"u}rel.
\newblock Relative {H}ofer-{Z}ehnder capacity and periodic orbits in twisted
  cotangent bundles.
\newblock {\em Duke Math. J.}, 123(1):1--47, 2004.

\bibitem{MR2324797}
Viktor~L. Ginzburg and Ba{\c{s}}ak~Z. G{\"u}rel.
\newblock The generalized {W}einstein-{M}oser theorem.
\newblock {\em Electron. Res. Announc. Math. Sci.}, 14:20--29, 2007.

\bibitem{MR2036376}
Norimichi Hirano and Slawomir Rybicki.
\newblock Some remarks on degree theory for {SO(2)}-equivariant transversal
  maps.
\newblock {\em Topol. Methods Nonlinear Anal.}, 22(2):253--272, 2003.

\bibitem{MR1984999}
Jorge Ize and Alfonso Vignoli.
\newblock {\em Equivariant degree theory}, volume~8 of {\em de Gruyter Series
  in Nonlinear Analysis and Applications}.
\newblock Walter de Gruyter \& Co., Berlin, 2003.

\bibitem{MR1330918}
W.~Klingenberg.
\newblock {\em Riemannian geometry}, volume~1 of {\em de Gruyter Studies in
  Mathematics}.
\newblock Walter de Gruyter \& Co., Berlin, second edition, 1995.

\bibitem{MR0478069}
Wilhelm Klingenberg.
\newblock {\em Lectures on closed geodesics}.
\newblock Springer-Verlag, Berlin, 1978.
\newblock Grundlehren der Mathematischen Wissenschaften, Vol. 230.

\bibitem{MR736839}
M.~A. Krasnosel{\cprime}skii and P.~P. Zabreiko.
\newblock {\em Geometrical methods of nonlinear analysis}, volume 263 of {\em
  Grundlehren der Mathematischen Wissenschaften [Fundamental Principles of
  Mathematical Sciences]}.
\newblock Springer-Verlag, Berlin, 1984.
\newblock Translated from the Russian by Christian C. Fenske.

\bibitem{MR0226651}
John~W. Milnor.
\newblock {\em Topology from the differentiable viewpoint}.
\newblock Based on notes by David W. Weaver. The University Press of Virginia,
  Charlottesville, Va., 1965.

\bibitem{MR2250797}
Jos{\'e} Ant{\^o}nio~Gon{\c{c}}alves Miranda.
\newblock Generic properties for magnetic flows on surfaces.
\newblock {\em Nonlinearity}, 19(8):1849--1874, 2006.

\bibitem{MR0346846}
N.~A. Niki{\v{s}}in.
\newblock Fixed points of the diffeomorphisms of the two-sphere that preserve
  oriented area.
\newblock {\em Funkcional. Anal. i Prilo\v zen.}, 8(1):84--85, 1974.

\bibitem{MR730159}
S.~P. Novikov and I.~A. Taimanov.
\newblock Periodic extremals of multivalued or not everywhere positive
  functionals.
\newblock {\em Dokl. Akad. Nauk SSSR}, 274(1):26--28, 1984.

\bibitem{MR0322899}
Frank Quinn and Arthur Sard.
\newblock Hausdorff conullity of critical images of {F}redholm maps.
\newblock {\em Amer. J. Math.}, 94:1101--1110, 1972.

\bibitem{MR2208800}
Felix Schlenk.
\newblock Applications of {H}ofer's geometry to {H}amiltonian dynamics.
\newblock {\em Comment. Math. Helv.}, 81(1):105--121, 2006.

\bibitem{MR0353372}
Carl~P. Simon.
\newblock A bound for the fixed-point index of an area-preserving map with
  applications to mechanics.
\newblock {\em Invent. Math.}, 26:187--200, 1974.

\bibitem{MR0185604}
S.~Smale.
\newblock An infinite dimensional version of {S}ard's theorem.
\newblock {\em Amer. J. Math.}, 87:861--866, 1965.

\bibitem{MR1133303}
I.~A. Taimanov.
\newblock Non-self-intersecting closed extremals of multivalued or
  not-everywhere-positive functionals.
\newblock {\em Izv. Akad. Nauk SSSR Ser. Mat.}, 55(2):367--383, 1991.

\bibitem{MR1185286}
I.~A. Taimanov.
\newblock Closed extremals on two-dimensional manifolds.
\newblock {\em Uspekhi Mat. Nauk}, 47(2(284)):143--185, 223, 1992.

\bibitem{MR0464304}
A.~J. Tromba.
\newblock A general approach to {M}orse theory.
\newblock {\em J. Differential Geometry}, 12(1):47--85, 1977.

\bibitem{MR493919}
A.~J. Tromba.
\newblock The {E}uler characteristic of vector fields on {B}anach manifolds and
  a globalization of {L}eray-{S}chauder degree.
\newblock {\em Adv. in Math.}, 28(2):148--173, 1978.

\end{thebibliography}

\end{document}